\documentclass [12pt,oneside] {amsart}
\usepackage{latexsym,amsfonts,amssymb,amscd,fullpage}
\usepackage [all]{xy}
\usepackage{tikz}
\usetikzlibrary{fit}
\usetikzlibrary{arrows}
\usetikzlibrary{matrix}
\SelectTips{cm}{12}

\theoremstyle{plain}

\newtheorem* {Thm*} {Theorem}
\newtheorem* {Prop*} {Proposition}

\theoremstyle {definition}

\newenvironment{Pf}[1]{{\noindent\sc Proof #1:}}{\qed\\}
\newcommand {\specialmap} [4] {\text {$ #1\negmedspace : #2 #3 #4 $}}
\newcommand {\map} [3] {\specialmap {#1} {#2}{\to} {#3}}
\newcommand {\longmap} [3] {\specialmap {#1} {#2}{\longrightarrow} {#3}}
\newcommand {\isomap} [3] {\specialmap {#1} {#2}{\overset {\cong{\phantom{.}}} 
            {\longrightarrow}} {#3}}

\newcommand {\Hom} {\operatorname {Hom}}
\newcommand {\id} {\operatorname{Id}}

\newcommand {\at}[1] {\arrowvert_{#1}}

\renewcommand {\(} {\left(}
\renewcommand {\)} {\right)} 

\renewcommand {\geq} {\geqslant}

\newcommand {\sub} {\subseteq}

\newcommand {\CC} {\mathbb C}
\newcommand {\on}[1] {\operatorname{#1}}

\newcommand{\BB}{\mathbb{B}}

\renewcommand {\leq} {\leqslant}

\newcommand {\pt} {\on{pt}}

\newcommand {\Sn} {{S_n}}

\renewcommand {\SS} {\mathbb{S}}

\newcommand {\SVect} {{\on{SVect^f}}}
\newcommand {\Sym} {\on{Sym}}
\newcommand {\tisomap} [3] {\specialmap {#1} {#2}{\overset {\cong}
            {\Longrightarrow}} {#3}}

\newcommand {\tensor}{\otimes}

\newcommand {\ZZ} {\mathbb Z}

\newcommand {\inv}{^{-1}}

\newcommand {\mC} {\mathcal C}

\newcommand {\mV} {\mathcal V}
\newcommand {\mW} {\mathcal W}
\newcommand {\mO} {\mathcal O}
\newcommand {\mF} {\mathcal F}
\newcommand {\mG} {\mathcal G}
\newcommand {\mA} {\mathcal A}
\newcommand {\mB} {\mathcal B}

\newcommand {\mD} {\mathcal D}
\newcommand {\<} {\langle}
\newcommand {\kk} {\bold{k}}

\newcommand {\hbtimes}{\, {\widehat{\boxtimes}}\,}

\newcommand {\ul}[1]{\underline{#1}}

\newcommand {\sgn} {\on{sgn}}

\renewcommand{\theta}{\vartheta}
\newcommand {\ttr} {\mathbb Tr}

\newcommand {\lra} {\longrightarrow}
\renewcommand {\ul}[1] {\underline{#1}}
\def\Ac{\mathcal{A}}
\def\Bc{\mathcal{B}}
\def\Cc{\mathcal{C}}
\def\Kc{\mathcal{K}}

\def\Ec{\mathcal{E}}

\def\Pc{\mathcal{P}}

\def\Tc{\mathcal{T}}

\def\Vc{\mathcal{V}}

\def\1{{\bf 1}}
\def\lra{\longrightarrow}

\def\st{\operatorname{st}}

\title{ Symmetric and exterior powers of categories} 
\author{Nora Ganter}
\author{Mikhail Kapranov}
\thanks{Ganter (University Melbourne) was supported by NSF grant DMS-0504539. During part of
  the project, she worked at Colby College, Maine.\\
Kapranov (Yale University) was supported by NSF Grant DMS-0801198.
}
\date {\today}
\begin{document}
\maketitle
%\tableofcontents
%
%
%

\begin{abstract}
  We define symmetric and exterior powers of categories, fitting into
  categorified Koszul complexes. We discuss examples and
  calculate the effect of these power operations on the categorical
  characters of matrix 2-representations. 
\end{abstract}

\section*{Introduction}

The classical formalism of symmetric and exterior powers in linear algebra
plays a fundamental role in many areas of mathematics and physics:
from geometry of Grassmann manifolds to supersymmetry and second
quantization. In this paper we develop a ``categorification" of this formalism.
That is, we define symmetric and exterior powers not for vector spaces over
a given field $\kk$, but for $\kk$-linear {\em categories}.

\vskip .2cm

Tensor products of $\kk$-linear categories have been defined in several
contexts \cite{Kapranov:Voevodsky, Deligne:Tannakian,  Lyubashenko}.
They all implement the following basic desideratum. Let $A_1, A_2$ be
two associative $\kk$-algebras, and $\mathcal V_i = A_i\on{-mod}$
be the category of left $A_i$-modules. Then the ``categorical tensor product"
$\mathcal V_1 \boxtimes\mathcal V_2$ should have something to do with
the category $(A_1\otimes_{\kk} A_2)\on{-mod}$. Similarly for sheaves on
spaces etc. 

\vskip .2cm

Given a workable concept of $\boxtimes$, the definition of the $n$th symmetric
power of a category $\mathcal V$ is rather straightforward: this is the category
of $S_n$-invariant objects in $\mathcal V^{\boxtimes n}$. Here $S_n$ is
the symmetric group on $n$ letters.  In the physical language, this
corresponds to considering ``orbifold models" (with respect to $S_n$). 
In particular, let $X$ be a smooth projective variety over $\CC$ and
$\mathcal V$ be the category of coherent sheaves on $X$.
 Taking the Grothendieck groups of the symmetric powers,
 we get the space
\[
\mathcal F\,\,=\,\,\bigoplus_{n\geq 0} K(\Sym^n(\mathcal V)) \,\,\,=\,\,\,
\bigoplus_{n\geq 0} K^{S_n}(X^n)
\]
which has the flavor of the loop Fock space. In fact, if we understand $K$ as
topological K-theory, then it is an observation of Grojnowski
\cite{Grojnowski} that $\mathcal F$ is the irreducible representation of
the loop Heisenberg algebra corresponding to $H^*(X, \ZZ)$. 
This   extends the results of \cite{Grojnowski, Nakajima}
on the cohomology of the Hilbert schemes of a projective surface. 

\vskip .2cm

The categorical analog of exterior
powers and the corresponding concept of  ``anti-equi\-va\-ri\-ance" with
respect to $S_n$, are less obvious. The key issue is what should play the
role of the sign character $\sgn: S_n\to \{\pm 1\}$.  There are (at least)
two answers to this. 

\newtheorem*{naive-answer}{Naive answer}
\begin{naive-answer}
 The categorical analog of   $\sgn$ should be an element
of the second cohomology $ H^2(S_n, \kk^*)$, a group which was 
essentially described by
I. Schur \cite{Schur} in his studies of projective representations of $S_n$.
In fact, Schur's results imply that for $n\geq 4$ and $char(\kk)\neq 2$
this group contains exactly one nontrivial element $c$. 
\end{naive-answer}
This point of view is motivated by the fact that the characters of any group
$G$ can be viewed as first homology classes: $\Hom(G, \kk^*)=H^1(G, \kk^*)$. 
  In physical language,  it corresponds to considering
``orbifold models with discrete torsion" \cite{Dijkgraaf}. 
Unfortunately,  literal implementation of this point of view leads to a formalism which
lacks some of the beauty and flexibility of the classical ``$S$-$\Lambda$-duality"
(for example, Koszul complexes). In this paper we adopt the following

\newtheorem*{super-answer}{Super-algebra answer}
\begin{super-answer}
The categorical analog of $\sgn$ should be a ``Picard character" (see Definition 3.1.3)
$$
  \on{Sgn}: S_n \to \Pc ic^{\ZZ/2}(\kk)
$$ 
with values in 
the category of super-lines over $\kk$.
Here by a super-line we mean a super-vector space of super-dimension
$(1|0)$ or $(0|1)$. Such a datum should combine both the classical character
$\sgn$ and the 2-cocycle $c$ as above. 
\end{super-answer}
This point of view draws on the wisdom accumulated in the theory of projective
representations of $S_n$ during the last 30 years \cite{Jozefiak:Projective}:
the theory becomes much more transparent if one introduces super-objects. 
Modern textbooks start with the super-approach right away \cite{Kleshchev}.

Using this as a guiding principle for our  definition of categorical exterior powers,
we find that several features of the classical theory generalize to the categorical
setting. In particular, we have categorical analogs of Koszul complexes
(\S 4) and we prove (Theorem 4.1.2)  that on the level of complexified
 Grothendieck groups they
give exact complexes of vector spaces. 

\vskip .2cm

The effect of categorical symmetric and exterior powers on Grothendeck groups
does not at all reduce to taking ordinary symmetric and exterior powers
of vector spaces. Rather, the generating functions for the dimensions of
the Grothendieck groups of these powers 
lead to expressions involving the Euler function $\phi(q)=\prod_{n\geq 1}(1-q^n)$
and remindful of the elliptic genus in topology. The reciprocity between
such generating functions given by categorical Koszul complexes 
includes, in the simplest instance, the relation between 
$\phi(q)^{-1}$, the generating function for the numbers of partitions and
$\prod_{n>0} (1+q^n)$, the generating function for the numbers of
strict partitions of integers. In more sophisticated representation-theoretic examples
(\S 5) the effect of categorical symmetric vs. exterior powers on
Grothendieck groups can be expressed via the ``logical proportion"
\[
\frac{\text{Symmetric powers}}{\text{Exterior powers}}\,\,\sim\,\,
\frac{\text{Untwisted Kac-Moody algebras}}{\text{Twisted Kac-Moody algebras}}.
\]

\vskip .2cm

There is a conceptual reason for the appearance of super-objects in 
projective representations of symmetric groups which, as far as we know,
has not been recognized in representation-theoretical
literature.  It comes from the
structure of the spherical spectrum $\SS$ in homotopy theory and from the
Barratt-Priddy-Quillen (BPQ)  theorem \cite{Priddy, Barratt:Priddy} which reconstructs this spectrum from
the symmetric groups.
 More precisely, we have the third (and truly fundamental)
answer to the above question ``what is the higher analog of the sign character":
\newtheorem*{third-answer}{Homotopy-theoretic answer}
\begin{third-answer}
There is a canonical homotopy theoretic character (coming from the BPQ theorem)
\[
\mathfrak{sgn}: S_n\lra \Omega\SS_0
\]
with values in the loop space of (one connected component of) the spherical spectrum. 
 The truncation of this loop
space in homotopy degrees 0 and 1 is described by the Picard category $\Pc ic^{\ZZ/2}(\ZZ)$
formed by super-lines with integer structure. The Picard character $\on{Sgn}$ above
is obtained from
$\mathfrak{sgn}$ by this truncation:
\[
\xymatrix{
S_n\ar[r]^{\mathfrak{sgn} }& \Omega\SS_0\ar@{~>}[r] &\Pc ic^{\ZZ/2}(\ZZ)\ar[r] &\Pc ic^{\ZZ/2}(\kk).
}
\]
\end{third-answer}
We explain this point of view in more detail in \S 3.1.9. As one can take less drastic truncations of
$\Omega\SS_0$, this provides a systematic way of defining  2-categorical, 3-categorical etc. analogs
of the sign character. For instance one could  take the truncation of $\Omega\SS_0$ in degrees
$[0,2]$  and interpret the resulting homotopy object in terms of a Picard 2-category,
thus getting a ``2-categorical sign character" from first principles. 

\vskip .2cm

The paper is organized as follows. In \S 1 we discuss various approaches to defining
tensor products of $\kk$-linear categories: abelian, pre-triangulated etc.  With the eye
on later applications, we also introduce the concept of a super-linear category
and discuss Grothendieck groups in this context. Symmetric powers of categories
are defined and studied in \S 2. In \S 3 we introduce exterior powers, starting from
an extended discussion in (3.1) of what should play the role of the sign character.
We consider first the ``naive" approach and then correct it using the super point of view. 
In \S 4 we introduce the categorical Koszul complexes and prove that they lead to
exact complexes of complexified Grothendieck groups. Further, \S 5 is devoted
to examples of symmetric and exterior powers coming from representation theory.
We show how they give companion Kac-Moody-type objects, one ``untwisted",
the other ``twisted". In \S 6 we study symmetric and exterior powers on the
category of 2-representations 
\cite{Ganter:Kapranov}of a group $G$  and find the effect of these constructions
on 2-characters of such 2-representations. Finally, \S 7 is devoted to discussion
of further directions and open question.

\vfill\eject

\section*{1. Tensor products of categories}

\vskip .5cm

\subsection*{(1.1) Definition of tensor products and direct sums.}

Let $\kk$ be an algebraically closed field.
By a linear category we will always mean an additive $\kk$-linear category, i.e.,
a category $\mathcal V$ with finite direct sums, 
in which all the sets $\on{Hom}_{\mathcal V}(V,W)$
are made into $\kk$-vector spaces so that the composition of morpisms is
bilinear. 

Let ${\mathcal V}, {\mathcal W}$ be two linear categories. 
Following \cite{Kapranov:Voevodsky}, \S 5.18
(see also \cite{Bakalov:Kirillov}, Def. 1.1.15),  we  define their tensor
product $\mV\boxtimes \mW$ to have, as objects, the symbols
(formal direct sums of formal tensor products)
$$\bigoplus_{i=1}^n V_i\boxtimes W_i, \quad V_i\in \on{Ob}(\mV), W_i\in\on{Ob}(\mW).
\leqno (1.1.1) $$
Morphisms between such symbols are defined as follows. First, we consider the case
when the formal sums in (1.1.1) have only one summand. In this case we put
$$\on{Hom}_{\mV\boxtimes\mW}\bigl( V\boxtimes W, \, V'\boxtimes W') \quad = \quad
\on{Hom}_{\mV}(V, V') \otimes_\kk \on{Hom}_{\mW}(W, W').\leqno (1.1.2) $$
 Then, we define
$$\on{Hom}_{\mV\boxtimes\mW}\left( \bigoplus_{j=1}^m V_j \boxtimes W_j,
 \,\,\, \bigoplus_{i=1}^n
V'_i\boxtimes W'_i \right) \quad =
\quad  \left\|\on{Hom}_{\mV\boxtimes\mW}\bigl( V_j\boxtimes W_j, \,
 V'_i\boxtimes W'_i)\right\|_{j=1, ..., m}^{i=1, ..., n}
\leqno (1.1.3)  $$
to consist of $n$ by $m$ matrices with the $(i,j)$th matrix element being a morphism from the $j$th
formal summand of the source to the $i$th formal summand of the target.  As this definition
of morphisms with matrices mimics morphisms between direct sums, we get that (1.1.1)
is indeed a categorical direct sum of
the $V_i\boxtimes W_i$ in $\mV\boxtimes\mW$. Further, notice the following:

\newtheorem* {1.1.4}  {(1.1.4) Proposition}
\begin{1.1.4}  Let $V_1, V_2$ be objects of $\mV$. Then
$$(V_1\oplus V_2) \boxtimes W \quad =  \quad (V_1\boxtimes W) \oplus (V_2\boxtimes W),$$ 
i.e., the left-hand side is a categorical direct sum, in
$\mV\boxtimes\mW$, of the two summands in the right-hand side.  
Similarly for $V\boxtimes(W_1\oplus W_2)$. 
\end{1.1.4}

\begin{proof} In any  additive category $\mathcal C$, to say  that an object $X$ is a categorical direct sum
 of $X_1$ and $X_2$, is the same as to say that there are  morphisms in $\mathcal C$
$$i_\nu: X_\nu\to X, \quad p_\nu: X\to X_\nu,\quad \nu =1,2,$$
such that 
$$p_\nu i_\nu = \on{Id}_{X_\nu},\quad  p_\nu i_\mu=0, \, \nu\neq \mu, \quad
i_1p_1+i_2p_2 = \on{Id}_X.$$
Suppose we have such morphisms for ${\mathcal C}=\mV$,  $X_\nu = V_\nu$ and $X=V_1\oplus V_2$.
Then tensoring them with $\on{Id}_W$ in (1.2), we get the required
 morphisms for ${\mathcal C}=\mV\boxtimes\mW$
and for the objects claimed in the proposition. 
\end{proof}

We also use the notation $\mV\boxplus\mW$ for the direct sum (Cartesian product)
of two linear categories, i.e., 
$$\on{Ob}(\mV\boxplus\mW) = \on{Ob}(\mV)\times \on{Ob}(\mW), \leqno (1.1.5)$$
$$
\on{Hom}_{\mV\boxplus\mW}\bigl( (V,W), (V', W')\bigr) = \on{Hom}_{\mV}(V, V')\oplus \on{Hom}_{\mW}
(W,W'). $$

\newtheorem*{1.1.6}{(1.1.6) Proposition}
\begin{1.1.6} For any three linear categories $\mV_1, \mV_2, \mW$ we have a natural equivalence
$$(\mV_1\boxplus\mV_2)\boxtimes \mW \quad = \quad (\mV_1\boxtimes W) \boxplus (\mV_2\boxtimes W).$$
Similarly for $\mV\boxtimes(\mW_1\boxplus \mW_2)$.
\end{1.1.6}

\begin{proof} A functor from the left-hand side to the right-hand side is given by
$$\Phi: (V_1, V_2)\boxtimes W \quad \mapsto \quad (V_1\boxtimes W, V_2\boxtimes W).$$
A functor from the right-hand side to the left-hand side is given by
$$\Psi: (V_1\boxtimes W_1, \, V_2\boxtimes W_2) \quad\mapsto\quad 
\bigl( (V_1, 0)\boxtimes W_1\bigr)\,\, \oplus \,\, \bigl( (0, V_2) \boxtimes W_2\bigr).$$
The fact that $\Phi$ and $\Psi$ are quasi-inverse to each other follows easily from
Proposition 1.1.4. \end{proof}

\subsection*{(1.2) Examples} We now list some immediate examples.

\vskip .2cm

\noindent {\bf (1.2.1) Algebras.} For a $\kk$-algebra $A$ we denote 
by $A\on{-mod}$ the category  of left $A$-modules and by $A\on{-mod}^f$
the subcategory of finitely presented left modules. Similarly, we write
$\on{mod-}A$ and $\on{mod^f-}A$ for the corresponding right module categories.
For any two algebras $A$ and $B$ we have a fully faithful embedding
\begin{eqnarray*}
(1.2.1)(a)\phantom{XXxiXXXX}
(A\on{-mod}^f)\boxtimes (B\on{-mod}^f) &\longrightarrow& (A\otimes_\kk
B)\on{-mod}^f   \phantom{(1.2.1)(a)XXxiXXXX}\\ 
\bigoplus_{i=1}^n M_i\boxtimes N_i &\longmapsto& \bigoplus_{i=1}^n M_i\otimes_k N_i,
\end{eqnarray*}
and similarly for right-modules.
We also have an equivalence
\begin{eqnarray*}
(1.2.1)(b)\phantom{XXXXXXX}
(A\on{-mod}^f)\boxplus(B\on{-mod}^f) &\longrightarrow& (A\oplus
B)\on{-mod}^f
\phantom{(1.2.1)(b)XXXXXXX}\\ 
\quad (M, N) &\longmapsto& M\oplus N, 
\end{eqnarray*}
and the analogous equivalence for right-module categories.
If $A$ and $B$ are semisimple then the embedding in (1.2.1)(a) is an equivalence
of categories. 

\vskip .2cm

\noindent {\bf (1.2.2) 2-vector spaces.}
By a 2-vector space we will mean a semisimple abelian $\kk$-linear category
with finitely many isomorphism classes of simple objects, cf.
\cite{Kapranov:Voevodsky}. 
Thus, if $A$ is semisimple, then $A\on{-mod}^f$ is a 2-vector space.
We will denote by $[n] = \kk^n\on{-mod}$
the standard $n$-dimensional 2-vector space. Each 2-vector space $\mV$
is equivalent to $[n]$, where $n$ is the number of simple objects in $\mV$ up
to isomorphism. We write $n=\on{Dim}(\mV)$.  Thus, we have
$$[n]\boxtimes [m] = [nm], \quad [n]\boxplus [m] = [n+m].$$

\vskip .2cm

\noindent {\bf (1.2.3) Coherent sheaves.} Let $X$ be a $\kk$-scheme of finite
type. Denote by $\on{Coh}_X$ the category of coherent sheaves of 
$\mO_X$-modules. For any two schemes $X, Y$ as above
we have a fully faithful embedding
\begin{eqnarray*}
  (1.2.3)(a)\phantom{XXXXXXXXXX} 
\on{Coh}_X\boxtimes\on{Coh}_Y &\longrightarrow& \on{Coh}_{X\times Y}
   \phantom{(1.2.3)(a)XXXXXXXXXXXXXX}\\
\mF\boxtimes\mG &\longmapsto &p_X^*\mF \otimes p_Y^*\mG,
\end{eqnarray*}
where $p_X, p_Y$ are the projections of $X\times Y$ to $X$ and $Y$, respectively. We also
have an equivalence
$$\on{Coh}_X\boxplus\on{Coh}_Y \longrightarrow \on{Coh}_{X\sqcup Y},
\leqno (1.2.3)(b)$$ 
where $X\sqcup Y$ is the disjoint union.
If $X= \on{Spec}(A)$ and $Y=\on{Spec}(B)$ are affine, this reduces to (1.2.1).

\vskip .2cm

We denote by $D^b(X) = D^b(\on{Coh}_X)$ the bounded derived category of
coherent sheaves on $X$. By passing to the derived functors in (1.2.3),
we get a fully faithful embedding and an equivalence
\begin{eqnarray*}
(1.2.3)(c)\phantom{XXXXXXXXX}
D^b(X)\boxtimes D^b(Y)&\longrightarrow& 
D^b(X\times Y)
\phantom{(1.2.3)(c)XXXXXXXXXXX}\\
D^b(X)\boxplus D^b(Y)&\longrightarrow& D^b(X\sqcup Y).
\end{eqnarray*}

\vskip .2cm

\noindent {\bf (1.2.4) Characterization by a (2-)universal property.} 
Let $\mathcal{X}$ be another additive $\kk$-linear category.
A functor
$$\psi: \mV \times\mW \to \mathcal{X}$$
is  called bilinear, if it is linear in each of the arguments. 
This implies, as in Proposition 1.1.4, 
 that  it takes direct sums in each argument to direct sums. 
For example, the functor
\begin{eqnarray*}
\boxtimes: \mV\times\mW&\longrightarrow& \mV\boxtimes \mW\\
 (V, W) &\longmapsto& V\boxtimes W
\end{eqnarray*}
is bilinear. Moreover, it is the universal bilinear functor, which means
that for any $\mathcal{X}$, precomposition with $\boxtimes$ defines an
equivalence 
$$
  \on{\mF\!un_\kk}(\mV\boxtimes\mW,\mathcal
  X)\stackrel\simeq\longrightarrow
  \on{\mF\!un_{\kk-bil}}(\mV\times\mW,\mathcal X)
$$
(categories of (bi-)linear functors and linear natural transformations).
%
%and $\psi$ as above there is a linear
%functor $f: \mV\boxtimes \mW\to\mathcal{X}$, unique up to
%a unique isomorphism, such that $\psi = f\circ\phi$.
This universal property,  noticed in \cite{Bakalov:Kirillov}, 
 characterizes the category $\mV\boxtimes\mW$
uniquely up to  an  equivalence of categories, which is  unique up to
a unique isomorphism.

\subsection*{(1.3) Completed tensor products: abelian case}
If the categories $\mV, \mW$ in (1.1) are abelian (resp. triangulated),
then $\mV\boxtimes\mW$ is, in general, no longer abelian (resp. triangulated).
One would like to enlarge $\mV\boxtimes\mW$ to a bigger abelian
(resp. triangulated) category $\mV\widehat{\boxtimes} \mW$ which would 
include kernels/cokernels (resp. cones) of morphisms between objects
of the form (1.1.1). 

Assume $\mV, \mW$ abelian, and let $\mathcal{X}$ be another abelian $\kk$-linear
category. We can then speak about bilinear functors
$$\Psi: \mV\times\mW\to\mathcal{X},$$
which are bi-right-exact, i.e., carry cokernels in each of the variables into
cokernels in $\mathcal{X}$. 
 The abelian tensor product of Deligne \cite{Deligne:Tannakian} is
an abelian category $\mV\hbtimes\mW$ together with a bi-right-exact functor 
$$\boxtimes: \mV \times\mW \to \mV\hbtimes\mW.\leqno (1.3.1)$$ satisfying
the following (2-)universal property:
for any other abelian category $\mathcal{X}$, precomposition with
$\boxtimes$ defines an equivalence 
$$
  \on{\mF\!un_{\kk,r.e.}}(\mV\hbtimes\mW,\mathcal
  X)\stackrel\simeq\longrightarrow
  \on{\mF\!un_{\kk-bil,r.e.}}(\mV\times\mW,\mathcal X).
$$
Here $\on{\mF\!un_{\kk(-bil),r.e.}}$ denotes the category of
  (bi-)linear (bi-)right exact functors and linear 
natural transformations. These properties characterize
$\mV\hbtimes\mW$, if it exists, uniquely up to an  equivalence of
categories, unique up to a unique isomorphism. 

\vskip .2cm

\noindent {\bf (1.3.2) Examples.} (a) (\cite{Deligne:Tannakian}, Prop. 5.3) 
Let $A$ and $B$ be right-coherent $\kk$-algebras. Then
$$(\on{mod^f-}A)\hbtimes (\on{mod^f-}B) \simeq \on{mod^f-}(A\otimes_\kk B).$$

\vskip .2cm

\noindent (b) (\cite{Deligne:Tannakian}, Prop. 5.11) Let $A$ be a finite-dimensional
$\kk$-algebra and $\mV$ any abelian $\kk$-linear category. Then $(\on{mod^f-}A)\hbtimes \mV$
exists and is equivalent to the category of $\kk$-linear right $A$-modules in $\mV$, 
i.e., objects $V\in\mV$ together with a $k$-algebra homomorphism
$A^{op}\to\on{End}_{\mV}(V)$. 

\vskip .2cm

\noindent (c) \cite{Lyubashenko} Let $(X, \mathcal{S})$ be a stratified space,
so $\mathcal{S}$ is the set of strata, partially ordered by inclusion of the closures.
 A perversity function for $(X,\mathcal{S})$
is a monotone function $p: \mathcal{S}\to \ZZ$. To this data one associates
an abelian category $\on{Perv}(X, \mathcal{S}, p)$ of $p$-perverse sheaves
of $\kk$-vector spaces on $X$ relative to $\mathcal{S}$. If $(X', \mathcal{S}')$
is another stratified space and $p'$ a perversity function for it then
the product $X\times X'$ is stratified by the products of the strata,
who form the set ${\mathcal{S}}\times\mathcal{S'}$. It has the perversity
function $p\dotplus p'$ given by
$$(p\dotplus p')(S\times S') = p(S) + p'(S'), \quad S\in\mathcal{S}, \, S'\in
\mathcal{S}'.$$
In this situation, we have
$$\on{Perv}(X, \mathcal{S}, p)\hbtimes \on{Perv}(X', \mathcal{S}', p') = 
\on{Perv}(X\times X', \mathcal{S}\times \mathcal {S}', p\dotplus p').$$

\vskip .3cm

\noindent {\bf (1.4) Completed tensor products: pre-triangulated case.} 
Let $\mV, \mW$, and  $\mathcal{X}$ be $\kk$-linear triangulated
categories. Then one can speak about bilinear functors
$\Psi: \mV\times\mW\to\mathcal{X}$
which are bi-exact, i.e., take exact triangles in each of the variables
into exact triangles in $\mathcal{X}$. A natural approach to the definition
of the ``triangulated tensor product'' of $\mV$ and $\mW$ would be to look
for the target of the universal bi-exact functor. However, as pointed out in
\cite{Bondal:Larsen:Lunts}, it is difficult to prove existence in
any nontrivial case, as one does not really have control over bi-exact functors.
An alternative approach, proposed in \cite{Bondal:Larsen:Lunts}, 
works with {\em enhanced triangulated categories}, following
\cite{Bondal:Kapranov}. 

\newtheorem*{141}{(1.4.1) Definition}
\begin{141}
By a dg-category we mean a category $\mV$ enriched in the monoidal category of
complexes of $\kk$-vector spaces. Thus for any $V, V'\in\on{Ob}(\mV)$ we have a complex
$\on{Hom}^\bullet_{\mV}(V, V')$. We denote by $H^0(\mV)$ the $\kk$-linear
category with the same objects of $\mV$ and
$$\on{Hom}_{H^0(\mV)}(V, V') \quad = \quad H^0 \on{Hom}_{\mV}(V, V').$$
A dg-functor  between dg-categories is a functor preserving the
enrichment, i.e., inducing morphisms of Hom-complexes. A dg-functor $F: \mV\to\mV'$
is called a quasi-equivalence, if $H^\bullet(F): H^\bullet(\mV)\to
H^\bullet(\mV')$ is an equivalence of categories. 
\end{141} 

\vskip .2cm

\noindent {\bf (1.4.2) Example.} If $\mA$ is an abelian $\kk$-linear
category, then $C(\mA)$, the category of cochain complexes
over $\mA$, is a dg-category (see \cite[p.94]{Bondal:Kapranov}). The
category $H^0C(\mA)$ is the homotopy category of complexes. 

\vskip .2cm

  We write $C(\kk\on{-mod})$ for the dg-category of complexes
  of $\kk$-vector spaces.  
  If $\mV$ and $\mW$ are two dg-categories, then the dg-functors from
  $\mV$ to $\mW$ form themselves the 
  objects of a dg-category
  $\on{\mF\!un_{dg}}(\mV,\mW)$
  (see \cite[p.95]{Bondal:Kapranov}) and the opposite category 
  $\mV^{\on{op}}$ is again a dg-category.
  For a dg-category $\mV$, we define
  $$\on{mod-}\!\mV:=
  \on{\mF\!un_{dg}}(\mV^{\on{op}},C(\kk-mod)).\leqno (1.4.3)$$ 
  This is a dg-category with shift and cone functors, inherited from
  $C(k\on{-mod})$. 
  The original dg-category $\mV$ is realized as a full
  dg-subcategory of $\on{mod-}\!\mV$ via the Yoneda embedding.

\newtheorem*{144}{(1.4.4) Definition}
\begin{144}
  Let $\on{Pre-Tr}(\mV)$ be the smallest full dg-subcategory of
  $\on{mod-}\!\mV$ that contains $\mV$ and is closed under isomorphisms,
  direct sums, shifts and cones.
  Let $\on{Perf}(\mV)$ be the full dg-subcategory of
  $\on{mod-}\!\mV$ consisting of semi-free\footnote{An object $F$ of
    $\on{mod-}\!\mV$ is {\em semi-free} if it has a filtration
    $0=F_0\subset F_1\subset\dots=F$, such that $F_{i+1}/F_i$ is
    isomorphic to a direct sum of shifted free objects, see
    \cite[p.11]{Bondal:Larsen:Lunts}.} 
  dg-modules that are homotopy
  equivalent to a direct summand of an object of $\on{Pre-Tr}(\mV)$.
\end{144} 

Bondal and Kapranov have given a hands-on description
of the category $\on{Pre-Tr}(\mV)$: 

\vskip .2cm

\newtheorem*{145}{(1.4.5) Proposition}
\begin{145}
The category $\on{Pre-Tr}(\mV)$ is 
equivalent to the category of {\em one sided twisted
  complexes}\footnote{In \cite{Bondal:Kapranov}, 
the category of one sided twisted complexes is called
$\on{Pre-Tr^+}(\mV)$.} of 
\cite[Sec.4]{Bondal:Kapranov}. 
%
%Let 
%$$
%  C= \bigoplus_{i=1}^m V_i [r_i]
%$$
%with
%$$
%  V_i\in\on{Ob}(\mV), \, r_i\in\ZZ,
%$$
%and let
%$$q =\|q_{ij}\|_{i<j}, \quad 
%\map{q_{ij}}{V_i[r_i]}{V_j[r_j+1]}$$
%be a triangular matrix of morphisms, satisfying $dq+q^2=0$. 
%Then $(C,q)$ defines an object $F_{(C,q)}$ of $\on{Pre-Tr}(\mV)$ as
%follows: for each 
%$W\in\on{Ob}(\mV)$, the value $F_{C,q}(W)$ is the graded abelian group
%$\bigoplus \Hom_\mV(W,V_i)[r_i]$ provided with the differential $d+q$.
%Every object of $\on{Pre-Tr}(\mV)$ is isomorphic to one of the form
%$F_{(C,q)}$. 
\end{145}

\begin{proof}
  In [loc.cit], Bondal and Kapranov prove that the category of
  one-sided twisted complexes may be embedded as a full subcategory in
  $\on{mod-}\!\mV$ in such a way that it contains $\mV$ and is closed
  under direct sums, shifts and cones. It follows that this category
  contains $\on{Pre-Tr}(\mV)$. To see equality, one notes that every
  one-sided twisted complex can be obtained from objects in $\mV$ by
  taking successive cones of degree zero morphisms (see
  \cite[Prop.3.10]{Bondal:Larsen:Lunts}). 
\end{proof}

We denote
$$\on{Tr}(\mV) = H^0 \on{Pre-Tr}(\mV). \leqno (1.4.6)$$

\newtheorem*{147}{(1.4.7) Definition}
\begin{147}
A dg-category $\mV$ is called {\em pre-triangulated}, if the embedding 
$$\longmap{i_{\mV}}\mV{\on{Pre-Tr}(\mV)}$$ is
a quasi-equivalence.
We say that $\mV$ is {\em perfect}, if
$\mV\to\on{Perf}(\mV)$ is a quasi-equivalence.
\end{147}

Note that ``perfect'' implies ``pre-triangulated''. 
If $\mV$ is pre-triangulated, then $H^0(\mV)$ is naturally a triangulated
category. Given a  triangulated cateogory $\mD$, an enhancement of $\mD$
is a pre-triangulated category $\mV$ with an equivalence 
$$\epsilon: H^0(\mV) \to\mD$$ of triangulated categories. 

\vskip .2cm

\noindent {\bf (1.4.8) Examples.} (a) If $\mV$ is any dg-category, then
$\on{Pre-Tr}(\mV)$ is pre-triangulated, and hence $\on{Tr}(\mV)$ is triangulated. 

\vskip .1cm

(b) If $\mA$ is any $\kk$-linear abelian category, then
$C(\mA)$, the category of complexes over $\mA$, is pre-triangulated. 

\vskip .1cm

(c) Let $X$ be a $\kk$-scheme of finite type. Then the triangulated
category $D^b(X)$ has the 
following enhancement $I(X)$. Let $\mO_X\on{-mod}$ be the abelian category of
all sheaves of $\mO_X$-modules (not necessarily quasicoherent). Then
$I(X)$ is the full dg-subcategory in $C(\mO_X\on{-mod})$ consisting of
complexes of injective objects, which are bounded below and have
only finitely many  cohomology sheaves, all coherent. The dg-category
$I(X)$ is perfect.  
%\marginpar{Question 3}
It is well known that $H^0I(X)$, i.e., the homotopy category of
complexes as above, is equivalent to $D^b(X)$. 

\vskip .2cm

%\marginpar{Question 4}
The following is an equivalent reformulation of \cite{Bondal:Larsen:Lunts},
Def. 4.2.

\newtheorem*{1.4.9}{(1.4.9) Definition}
\begin{1.4.9}
 Let $\mV, \mW$ be perfect dg-categories.
Their completed tensor product is defined by
$$\mV\hbtimes \mW \quad = \quad \on{Perf}(\mV  \boxtimes \mW),$$
where $\mV\boxtimes \mW$ is defined as in (1.1), using tensor products of
complexes and composition with a sign:
$$
  (f\tensor g)\circ(h\tensor k) :=
  (-1)^{\on{deg}(g)\cdot\on{deg}(h)}(fh)\tensor(gk). 
$$
\end{1.4.9}
\vskip .2cm

\noindent {\bf (1.4.10) Example.} (\cite{Bondal:Larsen:Lunts}, Th. 5.5)
 Let $X, Y$ be smooth projective varieties over $\kk$. Then
$I(X)\hbtimes I(Y)$ is quasi-equivalent
to $I(X\times Y)$. 
\vskip .2cm

There is no obvious (2-)universal property characterizing completed tensor
products of perfect dg-categories.

\vskip .3cm

\noindent {\bf (1.5) Super-objects. 2-periodic case.} Let
$(\on{SVect^f}, \otimes)$ 
be the symmetric monoidal category of finite dimensional 
super-vector spaces over $\kk$, see 
\cite{Manin} \cite{Kleshchev}. Thus, objects of $\on{SVect^f}$ are
finite dimensional $\ZZ/2$-graded
vector spaces $V= V_{\bar{0}}\oplus V_{\bar{1}}$, morphisms are linear operators
preserving grading, and $\otimes$ is the usual tensor product,
$$
  \(V\tensor W\)_{\bar 0} = V_{\bar 0}\tensor W_{\bar 0}\oplus
  V_{\bar 1}\tensor W_{\bar 1},$$
$$
  \(V\tensor W\)_{\bar 1} = V_{\bar 1}\tensor W_{\bar 0}\oplus
  V_{\bar 0}\tensor W_{\bar 1},
$$
with the symmetry
isomorphism 
given by the Koszul sign rule, 
\begin{eqnarray*}
  V\otimes W&\longrightarrow& W\otimes V\\
  v\tensor w&\longmapsto& (-1)^{\on{deg}(v)\cdot\on{deg}(w)}w\tensor v.
\end{eqnarray*}
The shift
of $\ZZ/2$-grading will be denoted by $\Pi$, so
$$(\Pi V)_{\bar{0}} = V_{\bar{1}}, \quad (\Pi V)_{\bar{1}} = V_{\bar{0}}. \leqno (1.5.1)$$
%
%%%%%Let $\on{SVect}^f$ be the category of finite-dimensional super-vector spaces.
%
\vskip .2cm

\newtheorem*{152}{(1.5.2) Definition}
\begin{152}
By a super-linear category we will mean a $\kk$-linear category which is a  module
category over $(\on{SVect^f}, \otimes)$.
\end{152}

Thus, in a super-linear category $\mV$ we have the functor $\Pi$ given by
$$\Pi (V) = \Pi(k)\otimes V,\leqno (1.5.3)$$
and each $\on{Hom}_\mV(V, V')$ is naturally extended to a super-vector space
with 
$$\on{Hom}_{\bar{0}}(V, W)  = \on {Hom}_{\mV}(V, W), \quad
\on{Hom}_{\bar{1}}(V, W) =  \on{Hom}_{\mV}(V, \Pi W). \leqno (1.5.4)$$ 
We will denote by $\mV_\bullet$ the {\em graded extension} of $\mV$, i.e., 
the category with the same objects as $\mV$ and
the extended ($\ZZ/2$-graded) $\Hom$-sets above.
We will refer to $\mV$ as the {\em even part} of $\mV_\bullet$.
An irreducible object $V$ of a superlinear category $\mV$ is called
{\em self-associate} if $V$ is isomorphic to $\Pi V$.  
Otherwise, $V$ is called {\em absolutely irreducible}.
%Absolute irreducibility is a necessary and, as it turns out,
%sufficient condition for $\bold V$ to remain irreducible in the
%extended category $((\SVect)^G)_\bullet$.

\vskip .2cm
\noindent {\bf (1.5.5) Examples.} (a) Let $A=A_{\bar{0}}\oplus
A_{\bar{1}}$ be a finite dimensional
associative and unitary superalgebra, i.e., an associative and unitary
algebra in $\on{SVect^f}$. A right $A$-supermodule $M$ is called {\em
  finitely generated} if there is an even surjective morphism 
$$
  A^{\oplus l}\oplus\Pi(A)^{\oplus k} \twoheadrightarrow M
$$ 
of right $A$-supermodules, and $M$ is called {\em finitely presented} if it fits into
an exact sequence of even morphisms
$$
  A^{\oplus i}\oplus\Pi(A)^{\oplus j} \longrightarrow 
  A^{\oplus l}\oplus\Pi(A)^{\oplus k} \longrightarrow M\longrightarrow
  0.
$$ 
We write $\on{Smod^f-}A$ for the category of finitely presented right
$A$-supermodules and even morphisms. This is a superlinear category. 
If $A$ is finitely generated
%{\em right-coherent}, i.e., if every finitely
%generated right superideal of $A$ is also finitely presented, 
then 
$\on{Smod^f-}A$ is an abelian superlinear category.  
%\marginpar{Question 5}
\vskip .1cm

(b) Consider the Clifford superalgebras 
$$
  C_1 = \kk[\xi]/(\xi^2)
$$
and
$$
  C_2 = C_1 \otimes C_1 \cong
  \kk[\xi_1,\xi_2]/(\xi_1^2,\xi_2^2,\xi_1\xi_2-\xi_2\xi_1), 
$$
see \cite[(12.1.3), (12.2.4)]{Kleshchev}. In $\on{Smod^f-}C_1$ we
have exactly one isomorphism class of irreducible objects, namely that
of $C_1$ viewed as right-module over itself. In particular, $C_1$ is
self-associate.  
In $\on{Smod^f-}C_2$, we have a decomposition 
$$
  C_2\cong U\oplus \Pi U,
$$
where $U$ and $\Pi U$ are non-isomorphic irreducible objects, and any
irreducible object is isomorphic to one of these two. It follows that
we have an equivalence of superlinear categories
$$
  \on{Smod^f-}C_2\medspace\simeq\medspace\SVect.
$$
\vskip .1cm

(c) Let now $\mV$ be a superlinear 2-vectorspace with $2n$ isomorphism
classes of absolutely irreducible objects (these come in pairs $U$,
$\Pi U$) and $m$ isomorphism classes of self-associate irreducible
objects. Then $\mV$ is equivalent, as a superlinear category, to the
{\em ``standard super 2-vectorspace''} of these dimensions,
$$
  [n|m] := \(\SVect\)^{n}\boxplus \(\on{Smod^f-}C_1\)^m.
$$
\vskip .1cm

(d) Let $\mA$ be an abelian $\kk$-linear category. Let $C^{(2)}(\mA)$ be the
category of 2-periodic complexes over $\mA$, i.e., of complexes
$(V^\bullet, d)$ such that $V^i = V^{i+2}$ and $d_i = d_{i+2}$ for all $i$.
Such complexes can be equally well considered $\ZZ/2$-graded. 
Morphisms in $C^{(2)}(\mA)$ are also assumed to be 2-periodic. Then
$C^{(2)}(\mA)$ as well as the homotopy category $H^0C^{(2)}(\mA)$
are super-linear categories with $\Pi$ induced by the shift of grading of complexes.

\vskip .1cm

(e) Let $\mV$ be a 2-periodic dg-category, i.e., each complex
$\on{Hom}_\mV^\bullet(V, V')$ is 2-periodic. 
Then $\mV^{\on{op}}$ and the category
$$\on{mod^{(2)}-}\mV:= \on{\mF\!un_{dg}}\(\mV^{\on{op}},C^{(2)}(\kk\on{-mod})\)$$ 
of contraviariant functors into even periodic cochain complexes are
again 2-periodic. So is
$\on{Pre-Tr}^{(2)}(\mV)$, the smallest full subcategory of
$\on{mod^{(2)}-}\mV$ that contains $\mV$ and is closed under isomorphisms,
direct sums, shifts and cones.

For each object $C$ we have $C[2] = C$, so the even subcategories of 
$\on{Pre-Tr}^{(2)}(\mV)$ and $\on{Perf}^{(2)} (\mV)$ are super-linear
with 
$\Pi C = C[1]$. 
Note that this definition of $\Pi$ makes sense in
$\on{Pre-Tr^{(2)}}(\mV)$, but that $\Pi$ might 
not be well-defined in $\mV$, even if the latter is pre-triangulated. In the
following, it is understood that the term
`{\em 2-periodic pre-triangulated} (or {\em perfect}) {\em
  dg-category}' refers to a
category $\mV$ of this kind where the objects $C[1]$
exist in $\mV$. Hence we may speak of the even subcategory
of $\mV$ as a superlinear category with graded extension $\mV$ (rather
than having to pass to $\on{Pre-Tr^{(2)}}(\mV)$).
%
%\vskip .1cm
%\marginpar{Question 3}
%(f) Let $X$ be a scheme of finite type  over $\kk$. Then $D^{(2)}(X)$,
%the derived category of 2-periodic complexes of coherent sheaves on $X$, is
%a super-linear triangulated category. In the case when $X$ is a
%smooth quasi-projective variety, it has a 2-periodic enhancement $I^{(2)}(X)$
%consisting of 2-periodic complexes of injective $\mO_X$-modules with
%coherent cohomology. The category $I^{(2)}(X)$ is a perfect super-linear
%category. 

\vskip .3cm

\noindent {\bf (1.6) Superlinear functors and supernatural transformations.} 
Let $\mV$ and $\mW$ be superlinear categories. A superlinear functor
from $\mV$ to $\mW$ consists of a linear functor 
$F\negmedspace :\mV\to\mW$ together with a linear natural isomorphism
$\phi\negmedspace : \Pi F\cong F\Pi$, satisfying
$$(\phi\Pi)\circ(\Pi\phi) = \id_F.$$ 
%We will view $\Pi$ itself as a superlinear functor
%where $\phi$ is taken to be $-\id$.
Given two such pairs, $(F,\phi)$ and $(H,\psi)$, a 
{\em supernatural transformation} between them is a linear natural
transformation $\xi\negmedspace :F\Rightarrow H$ satisfying 
$$
  \psi\circ(\Pi\xi) = (\xi\Pi)\circ\phi.
$$
We will denote the set of supernatural transformations from $(F,\phi)$
to $(H,\psi)$ 
by $$\on{s\mathcal N\!at}((F,\phi),(H,\psi)),$$
occasionally dropping $\phi$ and $\psi$ from the notation. 
Further, we will write 
$$
  s\ttr(F,\phi) := \on{s\mathcal N\!at}((\id,\id),(F,\phi))
$$
and
$$
  \on{sCenter}(\mV):= s\ttr(\id,\id).
$$   
\newtheorem*{161}{(1.6.1) Examples}
\begin{161}
(a) We have
$$
  \on{sCenter}(\SVect) \cong \kk.
$$
This includes diagonally into
$$
  \on{\mathcal N\!at}(\id,\id) \cong \kk^2. 
$$

\vskip .1cm
(b)
Similarly, 
$$
  s\ttr(\id,-\id) \cong \kk,
$$
but now the inclusion into $\on{\mathcal N\!at}(\id,\id)$ is the skew
diagonal map $a\mapsto(a,-a)$. This difference is picked up by the
action of $\Pi$: in the first case we have $\Pi\xi\Pi = \xi$, and in
the second case we have $\Pi\xi\Pi=-\xi$.

\vskip .1cm
(c)
More generally, let $\mV$ be a superlinear 2-vectorspace, and let $n$
be the number of isomorphism classes, up to shift of grading, of
irreducible objects in $\mV$. 
Then, by Schur's Lemma,
$$
   s\ttr(\id,\id)
   \cong \kk^{n},
$$
while the $\kk$-dimension of 
$$
   s\ttr(\Pi,\id)
$$
counts the isomorphism classes of self-associate objects of $\mV$.
\end{161}

\newtheorem*{162}{(1.6.2) Definition}
\begin{162}
  We will write $\on{s\mF\!un_\kk}(\mV,\mW)$ for the category of
  superlinear functors from $\mV$ to $\mW$ and supernatural
  transformations between them. This is itself a superlinear category
  with shift functor 
  $$\Pi(H,\phi):=(\Pi H,-\Pi\phi).$$
\end{162}
We have a superlinear equivalence of categories
$$
  \begin{array}{rcl}
  \SVect&\longrightarrow& \on{s\mF\!un_\kk}(\SVect,\SVect)
  \\
    V_\bullet&\longmapsto& (V_\bullet\tensor-,\phi_{V_\bullet}),
  \end{array}
  \leqno{(1.6.3)}
$$
where $\phi_{V_\bullet}$ is given by the symmetry isomorphism (with Koszul
sign) $$\phi_{V\bullet}\negmedspace:\Pi k\tensor V_\bullet\cong
V_\bullet\tensor \Pi k.$$ 

Superlinear functors extend canonically to even $\kk$-linear functors between the
extended categories, and we have an equivalence of categories
$$
  \on{s\mF\!un_\kk}(\mV,\mW)\stackrel\simeq\longrightarrow
  \on{\mF\!un_{\kk,ev}}(\mV_\bullet,\mW_\bullet), 
  \leqno{(1.6.4)}
$$
where $\on{\mF\!un_{\kk,ev}}$ stands for the category of 
$\kk$-linear functors that preserve the degree of
morphisms and $\kk$-linear natural transformations consisting of only
even maps.
\begin{proof}
  Indeed, a supernatural transformation between
  two superlinear functors is the same thing as an even $\kk$-linear natural
  transformation between the extended functors. 
  Further, let $\map F{\mV_\bullet}{\mW_\bullet}$ be an even
  $\kk$-linear functor, and let
  ${\id^\flat}\negmedspace :\Pi\Rightarrow \id$ be the natural isomorphism
  consisting of the identities, viewed as odd morphisms $\Pi V\to
  V$. Then $$F(\id^\flat)\negmedspace: F\Pi\Longrightarrow F$$ is an odd natural
  transformation, and we set $\phi:= F(\id^\flat)$, viewed as an even natural
  transformation from $F\Pi$ to $\Pi F$. Then the pair $(F\at\mV,\phi)$ is
  a superlinear functor, whose graded extension equals $F$.
\end{proof}
The composite of two super-linear functors is again a super-linear
functor in a canonical way, and composition is compatible with forming
the extended functors. 

\vskip .2cm
\noindent {\bf (1.6.5) Tensor products of superlinear categories:
uncompleted case.}
Let $\mV$ and $\mW$ be super-linear categories. 
Then the categorical tensor product of their extended categories is canonically
enriched over $\SVect$. 
Let 
$$
  \mV\boxtimes_s\mW := \(\mV_\bullet\boxtimes\mW_\bullet\)_{ev}.
$$
Note that the inclusion
$$
  \longmap B{\mV\boxtimes\mW}{\mV\boxtimes_s\mW}
$$
is faithful and essentially surjective, but not full. In particular,
in $\mV\boxtimes_s\mW$ we have the functor isomorphism
$$
  \id^\flat\boxtimes\medspace
  \Pi(\id^\flat):\Pi_{\mV_\bullet}\boxtimes
  \id_{\mW_\bullet}\Longrightarrow\id_{\mV_\bullet}\boxtimes
  \medspace\Pi_{\mW_\bullet},  
$$
which does not exist in $\mV\boxtimes\mW$. We make 
$\mV\boxtimes_s\mW$ into a superlinear category with shift functor
$\Pi\boxtimes\id$. It follows that the composite 
$$
  \longmap{B\boxtimes}{\mV\times\mW}{\mV\boxtimes_s\mW}
$$
is superlinear in both variables. The graded extension of
$\(\mV\boxtimes_s\mW,\Pi\boxtimes\id\)$ is naturally identified with
$\mV_\bullet\boxtimes\mW_\bullet$. Hence (1.6.4)
implies that $\boxtimes_s:=B\boxtimes$ is the universal such {\em
  bi-superlinear} functor out of $\mV\times\mW$. More precisely, the
pair $\(\mV\boxtimes_s\mW,\boxtimes_s\)$ 
satisfies the universal property (1.2.4) with (bi-)superlinear functors
and supernatural transformations in the place of (bi-)linear functors
and linear natural transformations.

\vskip .2cm
\noindent
\noindent {\bf (1.6.6) The abelian case.}
Let $\mA$ and $\mB$ be superlinear abelian categories. Then we may
apply the main Theorem of \cite{Greenough} (with $\mC=\SVect$) to
obtain the {\em `categorified coequalizer'} in the diagram

\begin{center}                                                                   
\begin{tikzpicture}                                                              
\matrix(m)[matrix of math nodes, row sep=-0.5cm, column sep=4em, text height=1.5ex, text depth=0.25ex]
{ {\phantom{\mA\hbtimes\mB}} & {\phantom{\mA\hbtimes\mB}} & {}\\
{\mA\hbtimes\mB} & {\mA\hbtimes\mB} & {\mA\medspace\!\widehat\boxtimes_s\medspace\!\mB}\\
{\phantom{\mA\hbtimes\mB}} & {\phantom{\mA\hbtimes\mB}} & {}\\};
\draw[->,font=\scriptsize,>=angle 90] (m-1-1) -- node[auto]                      
{$\Pi\widehat\boxtimes\mathrm{id}$} (m-1-2);                                                
\draw[->,font=\scriptsize,>=angle 90] (m-3-1) -- node[below]                     
{$\mathrm{id}\widehat\boxtimes\Pi$} (m-3-2);                                     
\draw[->,font=\scriptsize,>=angle 90] (m-2-2) -- node[auto] {$B$} (m-2-3);       
\end{tikzpicture}                                                                
\end{center} 

\noindent see \cite[Rem.3.6]{Greenough}. More precisely, a functor $F$ out of
$\mA\widehat\boxtimes\mB$ is called {\em SVect-balanced} if there is
a functor isomorphism 
$$
  \beta\negmedspace : F\circ(\Pi\boxtimes\id) \Longrightarrow F\circ
  (\id\boxtimes\medspace\Pi). 
$$
(This $\beta$ is part of the data of balanced functor, but
suppressed from the notation.)
The functor $B$ above is the universal right-exact SVect-balanced
functor out of $\mA\hbtimes\mB$. As in the previous
paragraph, we set
$\boxtimes_s:= B\boxtimes$, where
$\map{\boxtimes}{\mA\times\mB}{\mA\hbtimes\mB}$ is the
universal bi-right-exact, bilinear functor in (1.3.1). 
Since {\em bi-superlinear}
implies {\em SVect-balanced}, it follows that the pair
$({\mA\medspace\!\widehat\boxtimes_s\medspace\!\mB},\boxtimes_s)$ satisfies the
universal property in (1.3.1) with (bi-)superlinear right-exact
functors and supernatural transformations in the place of (bi-)linear
right-exact functors and linear natural transformations. 

%%
%%
%%  THE ARGUMENT IN THE OMITTED PROOF NEEDS TO INCLUDE A KOSZUL SIGN.
%% 
%%
%%
%%
\newtheorem*{167}{(1.6.7) Examples} 
\begin{167}
  (a) Let $A$ and $B$ be two 
  finitely generated
%right-coherent 
  superalgebras.
%  whose super-tensor product $A\otimes_\kk B$ is again right-coherent.
  The proof of \cite[Prop.5.3]{Deligne:Tannakian} can be formulated in
  the super setting (careful with Koszul signs!) to show that
%  let $\mC$ be
%  a superlinear abelian category, and let 
%  $$
%    \longmap{F}{\on{Smod^f-}A}{\mC}
%  $$
%  be a right-exact superlinear functor. Then 
%  $F(A)$ is a $\kk$-linear left $A$-supermodule in $\mC$, i.e.,
%  we have an even map of superalgebras 
%  $$
%    A\cong\on{End_\bullet}(A)\longrightarrow\on{End_\bullet}(F(A)),
%  $$
%  (endomorphisms of $A$ as right $A$-supermodule). 
%  In fact, we have an equivalence of categories
%  \begin{eqnarray*}
%    \on{s\mF\!
%      un_{r.e.}}(\on{Smod^f-}A,\mC)&\stackrel\sim\longrightarrow &
%    A-\mC\\ 
%    F&\longmapsto &F(A),
%  \end{eqnarray*}
%  where $A-\mC$ is the category of left $A$-supermodules in $\mC$ and
%  their morphisms. 
%
%  Similarly, if $B$ is a second right-coherent superalgebra, and we
%  are given a bi-superlinear, bi-right-exact functor
%  $$
%    \longmap{F}{\(\on{Smod^f-}A\)\times\(\on{Smod^f-}B\)}{\mC}
%  $$ 
%  then $F(A,B)$ carries
%  compatible structures as left $A$- and $B$-supermodule in $\mC$, and
%  we have an
%  an equivalence of categories 
%  \begin{eqnarray*}
%    \on{s\mF\!
%      un_{r.e.}^{bil}}\(\(\on{Smod^f-}A\)\times\(\on{Smod^f-}B\),\mC\)
%    & \stackrel\sim\longrightarrow & (A\tensor B)-\mC\\
%    F&\longmapsto& F(A,B).
%  \end{eqnarray*}
%  If $A\tensor B$ is again coherent, it follows that we have
%  $$
%    \(\on{Smod^f-}A\)\widehat\boxtimes_s\(\on{Smod^f-}B\)
%    \simeq
%    \on{Smod^f-}(A\tensor B),
%  $$
%  with the universal bi-superlinear bi-right-exact functor
%  $\boxtimes_s$ equal to
  \begin{eqnarray*}
    \(\on{Smod^f-}A\)\times\(\on{Smod^f-}B\)&\longrightarrow&
    \on{Smod^f-}\(A\tensor B\)\\
    (M,N) &\longmapsto& M\tensor_\kk N.
  \end{eqnarray*} 
makes $\on{Smod^f-}(A\tensor B)$ the categorical tensor product of
$\on{Smod^f-}A$ and $\on{Smod^f-}B$.

\vskip .1cm
\noindent
(b) Consider the 2-vectorspaces
$$
  [0|1] = \on{Smod^f-}C_1
\quad\quad\text{and}\quad\quad
  [1|0] = \SVect \simeq \on{Smod^f-}C_2. 
$$
By the previous example, 
$$[0|1]\medspace\widehat\boxtimes_s\medspace [0|1]\simeq [1|0].$$

\vskip .1cm
\noindent
(c)
Let $\mV$ be any superlinear abelian category. Then 
$$
  [1|0]\medspace\boxtimes_s\medspace \mV \simeq \mV
$$
is already abelian, so we do not need to complete the tensor product.
We arrive at the following general formula for completed tensor products
of superlinear 2-vectorspaces:
$$
  [k|l]\medspace\widehat\boxtimes_s\medspace [m|n] \quad\simeq\quad
  [km+ln|kn+lm].
$$

\vskip .1cm
\noindent
(d) Consider the $n$th Clifford superalgebra $C_n$ 
defined  by odd generators $\xi_1, ..., \xi_n$ subject to the relations
$$\xi_i^2 = 1, \quad \xi_i\xi_j  = -\xi_j\xi_i, \quad i\neq j.$$
Then we have $$C_n\cong C_1^{\tensor n}.$$
We recover the well known fact
$$
  \on{Smod^f-}C_n\simeq  \(\on{Smod^f-}C_1\)^{\widehat\boxtimes_s n}\simeq
  \begin{cases}
    [1|0] & \text{if $n$ is even, and}\\
    [0|1] & \text{if $n$ is odd.}
  \end{cases}
$$

\vskip .1cm
\noindent
(e) Let $\mV$ be a superlinear abelian category, and let $A$ be a
superalgebra of finite dimension over $\kk$. Then the argument in
\cite[5.11]{Deligne:Tannakian} is reformulated in the super-setting to
identify $\mV\medspace \!\widehat\boxtimes_s\medspace\!\(\on{Smod^f-}A\)$ 
with the category of right $A$-supermodules inside $\mV$.
\end{167}

\vskip .2cm
\noindent {\bf (1.6.8) The pre-triangulated case.} 
Let $C^\bullet$ and $D^\bullet$ be two 2-periodic complexes over
$\kk$. Then their tensor product as super-vectorspaces,
endowed with the differentials
$$
  d^{\bar 0}=\(
  \begin{matrix}
    d^{\bar 0}\tensor\on{id} &-\on{id}\tensor \medspace d^{\bar 1}\\
    \on{id}\tensor \medspace d^{\bar 0}&d^{\bar 1}\tensor\on{id}
  \end{matrix}\)
  \quad\quad\text{and}\quad\quad
    d^{\bar 1}=\(
  \begin{matrix}
    d^{\bar 1}\tensor\on{id} &\on{id}\tensor \medspace d^{\bar 1}\\
    -\on{id}\tensor \medspace d^{\bar 0}&d^{\bar 0}\tensor\on{id}
  \end{matrix}\),
$$
is again a 2-periodic complex.
Let $\mV$ and $\mW$ be 2-periodic perfect dg-categories, and
view $\mV\boxtimes\mW$ as 2-periodic dg-category with
the differentials just defined and composition with a sign as in (1.4.9).
We define the completed supertensor product of $\mV$ and $\mW$ as
$$
  \mV\widehat\boxtimes_s\mW := \on{Perf^{(2)}}(\mV\boxtimes\mW).
$$
%
%\vskip .2cm
%\noindent {\bf (1.6.9) Example.} \marginpar{Question 7}
% Let $X, Y$ be smooth projective varieties over $\kk$. Then
%$$I^{(2)}(X)\medspace\widehat\boxtimes_s\medspace I^{(2)}(Y)$$ is
%quasi-equivalent to $I^{(2)}(X\times Y)$. 

\vskip .3cm

\noindent {\bf (1.7) Grothendieck groups.} If $\mV$ is an abelian category
the Grothendieck group $K(\mV)$ is defined, as usual, by generators $\langle V\rangle$,
$V\in \on{Ob}(\mV)$ subject to the relations $\< V'\rangle  + \<V''\rangle  = \<V\rangle$ for each
exact sequence
$$0\to V'\to V\to V''\to 0.$$
If $\mV$ is an abelian superlinear category, then we define $K(\mV)=K^0(\mV)$
in a similar way, but imposing the additional relations
$$\<\Pi V\rangle   = - \<V\rangle, \quad V\in\on{Ob}(\mV). \leqno (1.7.1)$$
%\marginpar{Question 13}
This differs from the conventions in
\cite{Jozefiak:Projective} and \cite[(12.18)ff]{Kleshchev} by a sign.
For $n\geq 1$, we let 
$$
  K^n(\mV) := K^0\(\mV\medspace\widehat\boxtimes_s\medspace
  \on{Smod^f-}C_n\), 
$$
where $C_n$ is the $n$th Clifford superalgebra. By (1.6.7)(d), this
defines an even periodic theory. 
%
%If $\mV$ is a symmetric \marginpar{do
%we need ``symmetric in a super sense?'' What about the example of
%$C_1$-modules?} 
%monoidal category, then the tensor product of modules makes 
%$$
%  K^\bullet(\mV) = K^{0}(\mV)\oplus K^{1}(\mV)
%$$
%into a supercommutative $\ZZ/2$-graded ring.
We denote by
$$K^\bullet_\CC(\mW) = K^\bullet(\mW)\otimes \CC\leqno (1.7.2)$$
the complexified Grothendieck groups.
Assuming that there exists an $i\in\kk$ with $i^2=-1$, we may
define a product map
$$
  \longmap{\circledast}{K^\bullet(\mV)\otimes
    K^\bullet(\mW)}{K^\bullet(\mV\medspace\!\widehat\boxtimes_s\medspace\!\mW)}
  \leqno (1.7.3)
$$
as in \cite[(12.21)]{Kleshchev}, as follows: 
if 
$\on{deg}\langle V\rangle\cdot\on{deg}\langle W\rangle = 0$, we set
$$
  \langle V\rangle\circledast\langle W\rangle :=\langle V\boxtimes_sW\rangle.
$$
If $\on{deg}\langle V\rangle) = \on{deg}\langle W\rangle = 1$, we use
(1.6.7) to view $V\boxtimes_s W$ as a right $C_2$-supermodule object
inside $\mV\medspace\!\widehat\boxtimes_s\medspace\!\mW$.
It follows that we have the idempotent endomorphism 
$$
  f:=\frac{1+i\xi_1\xi_2}2
$$
of $V\boxtimes_s W$ with 
$$
  \on{ker}(f)\xi_1 = \on{im}(f),
$$
and we set
$$
  \langle V\rangle\circledast\langle W\rangle :=\langle \on{im}(f)\rangle.
$$

\vskip .2cm

\noindent {\bf (1.7.4) Examples.} (a) Let  $A$ be an associative superalgebra
over $\kk$,  which is coherent on the right. Then 
the category $\on{Smod^f-}A$ of finitely presented right $A$-supermodules
is a superlinear abelian category. 
On the other hand, $A= A_{\bar{0}}\oplus
A_{\bar{1}}$ can be considered as an ordinary associative algebra. Assuming that
this algebra is coherent on the right as well, we have the $\kk$-linear
abelian category 
$\on{mod^f-}A$ of finitely presented right $A$-modules. 
The Grothendieck groups of $\on{Smod^f-}A$ and $\on{mod^f-}A$ may be
different. 

\vskip .1cm
\noindent
(b) An instructive example is provided by the Clifford superalgebras
themselves. 
Considered as just an associative algebra, $C_n$ gives
$$K(\on{mod}^f-C_n) \quad = \quad\begin{cases} \ZZ & \text{if $n$ is
    even, and}\\
\ZZ^2 & \text{if $n$ is odd.}\end{cases}$$
as $C_n$ is either simple or has two simple summands, depending on the
parity of $n$. However, considered as a superalgebra, $C_n$ is simple,
see \cite[Exa.12.1.3]{Kleshchev}, and we have 
$$
  K^\bullet(\on{Smod^f-}C_n) = 
  \begin{cases}
    \ZZ\oplus(\ZZ/2)&\text{if $n$ is even, and}\\
    (\ZZ/2)\oplus\ZZ&\text{if $n$ is odd.}
  \end{cases}
$$
So, the complexified $K$-groups
$K^\bullet_\CC(\on{Smod^f-}C_n)$ are one-dimensional supervectorspaces
over $\CC$.
Simplifications of this kind will be crucial later in the paper. 

\vskip .1cm
\noindent
(c) Let $A$ and $B$ be finite dimensional superalgebras. Then 
$$
  \longmap{\circledast}{K_\CC^\bullet(\on{Smod^f-}A)\otimes
    K_\CC^\bullet(\on{Smod^f-}B)}
   {K_\CC^\bullet(\on{Smod^f-}(A\tensor B))}
$$
is an isomorphism (compare \cite[Lemma 12.2.15]{Kleshchev}).

\vskip .2cm

If $\mD$ is a triangulated $\kk$-linear category, then one defines
$K(\mD)$ by generators $\<V\rangle$, $V\in\on{Ob}(\mV)$,
and relations $\<V'\rangle +\< V''\rangle  = \<V\rangle$ for each
exact triangle
$$V'\to V\to V''\to V'[1].$$
Note that this implies automatically that
$$\<V[1]\rangle  = -\<V\rangle. \leqno (1.7.5)$$
%\marginpar{Question 8}
If $\mV$ is pre-triangulated, then we define
$$K(\mV) := K(H^0(\mV)).$$
Note that if $\mV$ is a 2-periodic pre-triangulated category, then
it is super-linear, with $\Pi V = V[1]$, so (1.7.4) can be
written as (1.7.1). 

\vskip .7cm

\vfill\eject 

\section*{2. Symmetric  powers.}

\vskip .5cm

\subsection *{(2.1) Reminder on 2-representations.} 
We keep working over the field $\kk$, as in (1.1).
Let $G$ be a group. A 2-representation (or categorical
representation) of $G$ is an action $\varrho$  of $G$ on a linear category $\mV$.
Explicitly, this means that we have the following data,
cf. \cite{Ganter:Kapranov}:  

\vskip .2cm

\noindent (2.1.1)(a) For each element $g\in G$, a linear functor $\varrho(g): \mV\to\mV$. 

\vskip .2cm

 \noindent (2.1.1)(b) For any pair of elements $(g,h)$ of $G$ an isomorphism of functors
      $$
        \tisomap{\phi_{g,h}}{(\varrho(g)\circ\varrho(h))}{\varrho(gh)}
      $$
 
\vskip .1cm

\noindent (2.1.1)(c)  An isomorphism of functors
      $$
        \tisomap{\phi_1}{\varrho(1)}{\id_c}
      $$
  \vskip .2cm

  such that the following conditions hold:

\vskip .2cm

\noindent (2.1.1)(d)   For any $g,h,k\in G$ we have
      $$
        \phi_{(gh,k)}(\phi_{g,h}\circ\varrho(k)) = \phi_{(g,hk)}(\varrho(g)\circ\phi_{h,k})
      $$
      (associativity); we also write $\phi_{g,h,k}$,

\vskip .2cm
   
\noindent (2.1.1)(e) We have
      $$
        \phi_{1,g} = \phi_1\circ\varrho(g) \quad\text{and}\quad
        \phi_{g,1} = \varrho(g)\circ\phi_1.
      $$

 \vskip .2cm
\noindent {\bf (2.1.2) Examples.} (a) Let 
$$c: G\times G\to \kk^*$$
 be a 2-cocycle,
i.e., a function satisfying the identity
$$c(g h, k)\cdot c(g, h) = c(g, hk)\cdot c(h, k).$$
We then have an action $\varrho=\varrho_c$ of $G$ on $\on{Vect}^f$.
By definition, for $g\in G$ the functor $\varrho(g): \on{Vect}^f \to \on{Vect}^f$
is the identity, and the  transformation
$\phi_{g,h}:  \id\Rightarrow \id$
is multiplication with $c(g,h)$, while $\phi_1$ is the
multiplication by $c(1,1)$.

The 2-representation $\varrho_c$
is equivalent to $\varrho_{c'}$ (in the appropriate sense, see 
\cite[2.2]{Bartlett}) 
if and only if $c$ is cohomologous to $c'$,
and $H^2(G, \kk^*)$ is identified with the set of ``1-dimensional 2-representations
of $G$'', i.e., $G$-actions on $\on{Vect}^f$ modulo equivalence. 

\vskip .2cm

(b) Alternatively, let 
$$
  1\longrightarrow \kk^*\longrightarrow\widetilde
  {G}\stackrel\pi\longrightarrow G\longrightarrow  1
$$
be a   central extension of $G$. For every
$g\in G$ the set $\pi\inv(g)$ is a $\kk^*$-torsor, and therefore
$$
  L_g := \pi\inv(g)\cup\{0\}
$$
is a $1$-dimensional $\kk$-vector space. Note that $L_1=\kk$. 
 The group structure on $\widetilde
{G}$ induces isomorphisms
$$
\mu_{g,h}:   L_g\tensor_k L_h\longrightarrow L_{gh}.
$$
satisfying the associativity conditions. For each $g\in G$ we then have
the functor of tensor multiplication with $L_g$:
\begin{eqnarray*}
  \varrho(g): \on{Vect}^f&\longrightarrow&\on{Vect}^f\\
  V & \longmapsto & L_g\otimes V.  
\end{eqnarray*}

The isomoprhisms $\mu_{g,h}$ give then the isomorhisms 
$\phi_{g,h}$ from (2.1.1)(b), and we define the functor $\phi_1$
to be the identity. This gives an action of $G$ on $\on{Vect}^f$.
One sees easily that if $\widetilde{G}$ corresponds to a cocycle $c$, then
this 2-representation of $G$ is equivalent to $\varrho_c$ from (a).  

\vskip .2cm

(c) A superlinear category is the same as a $\kk$-linear category with
a strict linear action by $\ZZ/2$. ``Strict'' means that all the
$\phi_{g,h}$ and $\phi_1$ are identity maps. The definitions of
{superlinear functor} and {supernatural transformation} in
(1.6) translate into the definition of 1- and 2-morphisms of
2-representations \cite[2.2.2]{Bartlett}. 
\vskip .2cm

\noindent {\bf (2.1.3)} Given a 2-representation of $G$ on $\mV$, we denote
by $\mV^G$ the {\em category of $G$-equivariant objects in $\mV$},
i.e., objects $V\in\mV$ equipped with isomorphisms
$$ \map{ \epsilon_g}{V}{\varrho(g)(V)},\quad g\in G,\leqno (2.1.4)$$
satisfying the following
compatibility condition: %first, it is required that for $g=1$ we have
%$$\epsilon_1 = \phi_{1,X}^{-1}: X\mapsto \varrho(1)(X).$$
%Second, it is required that 
for any $g,h\in G$ the diagram
$$\xymatrix{
V\ar[0,2]^{\epsilon_g}\ar[d]_{\epsilon_{gh}} &&
{\varrho(g)(V)} \ar[d]^{\varrho(g)(\epsilon_h)} \\
{\varrho(gh)(V)} &&
{\varrho(g)(\varrho(h)(V))}\ar[0,-2]^{\phi_{g,h,V}}
}$$
is commutative.

\vskip .2cm

\noindent {\bf (2.1.5) Example.}  Suppose that the $G$-action on $\mV$ is
trivial:
  all $\varrho(g)$ and $\phi_{g,h}$, as well as $\phi_1$,  are the identities.
Then a $G$-equivariant object in $\mathcal{V}$ is the same as a representation of $G$
in $\mathcal V$, i.e., an object $V\in\mathcal{V}$ and a homomorphism
$G\to \on{Aut}_{\mathcal {V}}(V)$. 

\vskip .2cm

\noindent {\bf (2.1.6)} If $\mV, \mW$ are two 2-representations of $G$, 
then $\mV\boxtimes\mW$ is again
a 2-representation, in an obvious way. Further, if $\mV, \mW$ are abelian
and $\mV \hbtimes \mW$ exists, then it is also a 2-representation of $G$.
This follows from the characterization of $\mV \hbtimes \mW$ by a (2-)universal
property. Similarly, if $\mV, \mW$ are perfect dg-categories with $G$-action,
then so is $\mV\hbtimes\mW$.

\vskip .2cm

\noindent {\bf (2.1.7) Example.} For 1-dimensional 2-representations
$\mV_c$ as in Example 2.1.2(a), we have
$$\mV_c \boxtimes \mV_{c'} \quad \simeq\quad  \mV_{c\cdot c'}.$$
So $H^2(G, \kk^*)$ is interpreted as the Picard group of 1-dimensional
2-representations, with operation given by $\boxtimes$. 

\vskip .3cm

\subsection*{(2.2) Symmetric powers.} Let $\mV$ be a linear category. We have then
the linear category
$$\mV^{\boxtimes n} = \mV \boxtimes ... \boxtimes \mV \quad (n \,\,\,\on{times}).$$
Let $S_n$ denote the symmetric group of permutations of $\{ 1, ..., n\}$.
Then we have the $S_n$-action on $\mV^{\boxtimes n}$ given by
$$\sigma(V_1\boxtimes ... \boxtimes V_n) \quad = \quad V_{\sigma\inv(1)} \boxtimes ... \boxtimes
V_{\sigma\inv(n)},$$
and extended to direct sums by additivity. 

\vskip .2cm

\newtheorem*{221}{(2.2.1) Definition}
\begin{221}
The $n$th symmetric power $\Sym^n(\mV)$ is the category of $S_n$-equivariant
objects in $\mV^{\boxtimes n}$:
$$\Sym^n(\mV) \quad =\quad (\mV^{\boxtimes n})^{S_n}.$$
\end{221}

\vskip .2cm

\noindent {\bf (2.2.2) Example.} Let $\mV$ be the standard 1-dimensional 2-vector
space, so $\mV = [1] = \on{Vect}^f$, see (1.2.2). Then $\mV^{\boxtimes n} = 
\on{Vect}^f$ and the $S_n$-action is trivial. So by Example 2.1.5, we have
that $$\Sym^n(\mV) = \on{Rep}(S_n)$$ 
is the category of finite-dimensional representations
of $S_n$ over $\kk$. Assume
$\on{char}(\kk) =0$.  Denote by $p(n)$ the number of partitions of $n$. Then, as well known,
$\on{Rep}(S_n)$ has $p(n)$ simple objects (irreducible representations of $S_n$).
So we can write
$$\Sym^n[1] \quad \simeq \quad [p(n)].\leqno (2.2.3)$$
This formula indicates that tensor operations on categories
lead to a new type of $\lambda$-rings, different from the
standard (special) $\lambda$-rings which, as well known,  satisfy
$$\Sym^n(1) = 1, \quad \forall n.$$
We will study such ``2-special'' $\lambda$-rings in a separate paper.

\vskip .2cm

\newtheorem*{224}{(2.2.4) Proposition}
\begin{224}
We have an equivalence
$$\Sym^n(\mV\boxplus \mW) \quad \simeq\quad \boxplus_{i+j=n} \,\, \Sym^i(\mV)\boxtimes
\Sym^j(\mW). $$
\end{224}

\noindent {\sl Proof:}  This follows from Proposition 1.1.6. \qed

\vskip .2cm

Let us denote by
$$\phi(q) =\( \sum_{n\geq 0} p(n) q^n\)\inv = \prod_{n\geq 1} (1-q^n)
\leqno (2.2.5)$$ 
the Euler function. 

\newtheorem*{226}{(2.2.6) Corollary}
\begin{226} Assume $\on{char}(\kk)=0$.
Let $\mV$ be a 2-vector space of dimension $d$. Then 
$$\sum_{n\geq 0} \on{Dim} \,\,\Sym^n(\mV) q^n \quad = \quad \frac1{\phi(q)^d.}$$
\end{226}

\vskip .2cm

\noindent {\bf (2.2.7)} If $\mV$ is abelian or perfect, then we can form the completed symmetric
product $\widehat{Sym}^n(\mV)$.  For this, we start with the completed
tensor power 
$$\mV^{\hbtimes n} \quad = \quad \mV \hbtimes \cdots \hbtimes \mV,$$
which we assume to exist in the abelian case, and define 
 as in (1.4.9) in the perfect
case. Then $\widehat{Sym}^n(\mV)$ is defined as the category of
$S_n$-equivariant objects in $\mV^{\hbtimes n}$. 

\vskip .2cm

%%%\marginpar{Qeustion 9}
\noindent {\bf (2.2.8) Examples.}(a)  Let $X$ be a smooth projective variety.
Then $\widehat{\Sym}^n (I(X))$ is the dg-category of
$S_n$-equivariant complexes of injective $\mO$-modules on $X^n$
that are bounded below and have  bounded coherent cohomology.
The corresponding $H^0$-category is the bounded derived category of
$S_n$-equivariant complexes of coherent sheaves on $X^n$. 

\vskip .1cm

(b) Let $A$ be a finite-dimensional $\kk$-algebra, and
$\mV=\on{mod^f-}A$. Then $V^{\widehat{\boxtimes} n}$ is the category
of right $A^{\otimes n}$-modules. The group $S_n$ acts on
$A^{\otimes n}$ by algebra automorphisms, so we have the crossed
product algebra $A^{\otimes n}[S_n]$, and
$$\widehat{Sym}^n(\mV) \quad = \quad \on{mod^f-}A^{\otimes n}[S_n].$$

\vfill\eject

\section*{3. Exterior powers.}

\subsection*{(3.1) Picard categories and categorical characters.}
 The construction of exterior powers in usual linear algebra is based on
two steps:
\begin{itemize}
\item[(1)] Identification of the sign character of $S_n$ as a homomorphism 
$\sgn: S_n\to\ZZ/2$ which is universal among homomorphisms $S_n\to A$
into abelian groups. At this stage $\ZZ/2$ is considered as an abstract 2-element group.

\item[(2)] Realization of $\ZZ/2$ as the subgroup $\{\pm 1\}\subset \kk^*$ of the multiplicative
group of the field, which is then used to multiply elements of various vector spaces.
\end{itemize}
We start by discussing the categorical analog of (1).  For this,  abelian groups $A$ 
should be replaced by Picard categories. We recall basic definitions and examples, following
\cite{SGA4} Exp. XVIII.

\newtheorem*{picard-cat-def}{(3.1.1) Definition}
\begin{picard-cat-def}
By a {\em (symmetric)  Picard category} we mean a 
(symmetric) monoidal category 
$(\mA, \otimes, \1)$ 
in which all
the objects are invertible, and all the morphisms are invertible with
respect to composition. 
\end{picard-cat-def}

 For any essentially small Picard category $\mathcal{A}$ we have
the  groups
$$\pi_0(\mathcal {A}) = ( \text{Ob}(\mathcal{A})/\text{iso}, \otimes), \quad \pi_1(\mathcal{A})=
\text{Hom}_{\mathcal{A}}(\mathbf{1}, \mathbf{1}).$$
The group $\pi_1(\Ac)$ is always abelian. If $\Ac$ is symmetric, then $\pi_0(\Ac)$ is also abelian. 

\newtheorem*{picard-cat-ex}{(3.1.2) Examples}
\begin{picard-cat-ex} (a)  Let $G$ be a  group. We can consider the set $G$ as
a discrete category (the set of objects is $G$, the only morphisms are identities). The
group operation in $G$ makes then $G$ into a Picard category, which we will
still denote
$G$ and call the {\em discrete Picard category} associated to $G$.
This Picard category has $\pi_0=G$, $\pi_1=0$.
If $G$ is abelian, then we get a symmetric Picard category. 

\vskip .2cm
 
 (b) Let $A$ be an abelian group.  Denote by  $A-\Tc \!ors$ 
 the category of $A$-torsors (i.e., principal homogeneous spaces $T$ over $A$)
 and their isomorphisms. This category is a symmetric Picard category with
 respect to the tensor product of torsors defined by
 $$
 T\otimes_A T' \,\,=\,\, (T\times T')\bigl/ \bigl\{(at, t') \sim (t, at'), \quad a\in A, \, t\in T, \, t'\in T'\bigr\}.
 $$
 Note that 
  $A-\Tc\! ors$ is equivalent to 
   the Picard category
 $\Bc A$ with one object $\1$ and
 $\Hom(\1, \1)=A$. The operation $\otimes$ on morphisms,
 as well as the composition of morphisms in $\Bc A$, are both given by the group operation in $A$. 
 Both $A-\Tc\! ors$ and $\Bc A$ have $\pi_0=0$, $\pi_1=A$. 
 
 \vskip .2cm
 
 (c) Let $A^\bullet = \{A^0\buildrel d\over\longrightarrow A^1\}$ be a two-term
complex of abelian groups. Then we have a symmetric Picard category $[A^\bullet]$ with
$$\text{Ob}([A^\bullet]) = A^1, \quad \text{Hom}_{[A^\bullet]}(a, a') = \{ b\in A^0:
d(b) = a-a'\}.$$
  Note that
$$
\pi_i[A^\bullet] = H^{1-i}(A^\bullet), \quad i=0,1. 
$$

\vskip .2cm

(d) Let $\kk$ be a field. We have then the symmetric Picard category $\Pc ic^\ZZ(k)$  of {\em graded lines}
over $\kk$. By definition, objects of  $\Pc ic^\ZZ(k)$
 are  $\ZZ$-graded $\kk$-vector spaces $V=\bigoplus_{n\in Z} V_n$
 of total dimension 1. In other words, such a $V$ has exactly one graded component
 of dimension 1, while all other graded components vanish. Morphisms in 
 $\Pc ic^\ZZ(k)$ are isomorphisms of graded vector spaces. The monoidal operation
 is given by the usual graded tensor product, while the symmetry is given by
 the Koszul sign rule:
 $$ x\otimes y \longmapsto (-1)^{\deg(x)\deg(y)} y\otimes x.
 $$
 We also denote by $\Pc ic^{\ZZ/2}(\kk)$ the  symmetric Picard category of {\em super-lines}
 over $\kk$ by which we mean
  super-vector spaces over $\kk$ of total dimension 1, i.e., of super-dimension
 either $(1|0)$ or $(0|1)$. We have
 \[
 \pi_0(\Pc ic^{\ZZ}(\kk))=\ZZ, \,\,\,\pi_1(\Pc ic^\ZZ(\kk))=\kk^*, \quad \pi_0(\Pc ic^{\ZZ/2}(\kk))=\ZZ/2, \,\,\,
 \pi_1(\Pc ic^{\ZZ/2}(\kk)) = \kk^*.
 \]

\end{picard-cat-ex}

\vskip .2cm

\newtheorem*{cat-char-def}{(3.1.3) Definition}
\begin{cat-char-def}
Let $G$ be a group and $\mA$ be a symmetric Picard category.  An $\Ac$-valued
{\em Picard character} of $G$ is a monoidal functor $X: G\to\mA$ 
(where $G$ is considered as a discrete Picard category). 
\end{cat-char-def}

Explicitly, a Picard character consists of objects $X(g)\in\Ac$ given for
each $g\in G$, of isomorphisms $$\phi_{g,h}: X(g)\otimes X(h)\to X(gh),$$ and
of an isomorphism $\phi_1: X(1)\longrightarrow \1$ which satisfy the conditions
identical to (2.1.1)(d-e), with $\circ$ replaced by $\otimes$. 

In fact, both Picard characters and 2-representations
as defined in (2.1),  are
particular cases of a more general concept: a 2-representation of $G$ in
a 2-category $\Cc$, see \cite{Ganter:Kapranov}   Def. 4.1. More precisely,
the context of (2.1) corresponds to $\Cc$ being the 2-category of $\kk$-linear
categories, their linear functors and natural transformations.
 Definition (3.1.3) corresponds to the case when 
$\Cc$ is the 2-category with one object $\pt$
canonically associated to the monoidal category $\Ac$: here
$$\on{1Hom_\Cc}(\pt, \pt) = \Ac$$ and composition of 1-morphisms corresponds
to the monoidal structure on $\Ac$.

\newtheorem*{cat-char-ex}{(3.1.4) Examples}
\begin{cat-char-ex} 
(a) If $\Ac=A$ is a discrete symmetric Picard category, then an $A$-valued
Picard character is simply a group homomorphism $G\to A$. More generally, for any
symmetric Picard category $\Ac$ and any Picard character $G\to \Ac$ we obtain a
homomorphism $G\to\pi_0(\Ac)$ with an abelian group target. 

\vskip .2cm

(b) If $\Ac=A-\Tc\! ors$ is the category of torsors over an abelian group $A$, then an
$\Ac$-valued Picard character of $G$ is the same as a central extension
$$
1\to A\lra \widetilde G\buildrel p\over \lra G\to 1.
$$
Indeed, given such an extension, every preimage $$X(g)=p^{-1}(g)$$ is
a torsor over 
$A=p^{-1}(1)$, and the group operation in $\widetilde G$ gives identifications
$$X(g)\otimes_A X(G')\longrightarrow X(gg').$$
Conversely, given any Picard character
$X: G\to \Ac-\Tc\! ors$, the set $$\widetilde G = \coprod_{g\in G}
X(g)$$ becomes a group, fitting into a central extension as
above. This generalizes Example 2.1.2(b) 
which corresponds to $A=\kk^*$. 

\vskip .2cm

(c) If $\Ac=\Bc A$  for an  abelian group $A$, then an $\Ac$-valued
Picard character is the same as a 2-cocycle $c: G\times G\to A$.

\end{cat-char-ex}

\noindent {\bf (3.1.5) Symmetric Picard categories and $[0,1]$-spectra.}
The construction of symmetric Picard categories in ExampleÊ(3.1.2)(c) is 
not the most general one. 
It is an insight of Grothendieck 
(see \cite{Drinfeld}  (5.5.2) for a discussion)
 that a complete description  
can be obtained in terms of  spectra, Êi.e., 
objects   of the stable homotopy category.\footnote{For an
  introduction to spectra, see \cite{Adams}.}  
 We recall some details of this description.
A {\em spectrum} $Y$ can be seen as an infinite loop space, i.e., a topological
space $\Omega^\infty Y$ together with a sequence of its deloopings
$\Omega^{\infty-j}Y$, $j\geq 0$,  which are topological spaces equipped  with
homotopy equivalences $\Omega^j(\Omega^{\infty-j}Y) \sim\Omega^\infty Y$. 
In particular, a spectrum $Y$ has homotopy groups  $\pi_i(Y)$, $i\in\ZZ$,
which are defined as
$$\pi_i(Y) \,\,=\,\,\pi_{i+j} (\Omega^{\infty -j}(Y)), \quad j\gg 0.
$$ 
By a $[0,1]$-{\em spectrum} we mean a spectrum $Y$
which has only two nonzero
homotopy groups, namely $\pi_0(Y)$ and $\pi_1(Y)$.

Note that any Picard category $\Ac$, considered as an abstract category, is a groupoid.
So for a small $\Ac$ the classifying space $B(\Ac)$ is a topological space with
$$\pi_0(B(\Ac))\,\,=\,\,\pi_0(\Ac), \quad \pi_1(B(\Ac), \{\1\}) = \pi_1(\Ac),$$
while the higher homotopy groups of each component of $B(\Ac)$ vanish. 

\newtheorem*{thm-picard-spectra}{(3.1.6) Theorem}
\begin{thm-picard-spectra}
 (a) The functor of taking the classifing space extends to a functor
$$
\BB: \bigl\{ \text{Small Symmetric Picard Categories} \bigr\} \lra \text{ \{$[0,1]$-spectra}\}
$$
such that $\pi_i(\BB(\Ac)) = \pi(\Ac)$, $i=0,1$. 

(b)  
The functor $\BB$ takes equivalences of
symmetric Picard categories into homotopy equivalence of $[0,1]$-spectra. 
After inverting the two types of equivalences, $B^\infty$ becomes an equivalence
of localized categories.  
 \end{thm-picard-spectra}
 
 For convenience of the reader we sketch the construction of the functor $\BB$.
 The infinite loop space of the spectrum $\BB(\Ac)$ is,  by definition,
 $\BB^\infty(\Ac)=B(\Ac)$, the usual classifying space. For any
 $m\geq 0$ the   $m$th delooping $\BB^{\infty -m}(\Ac)$ is defined  as
  the geometric
 realization of the simplicial set $\BB^{\infty -m}_\bullet(\Ac)$ constructed as follows. 
 By definition, $n$-simplices of $\BB^{\infty -m}_\bullet(\Ac)$ are
 ``$\Ac$-valued $m$-cocycles on the standard simplex $\Delta[n]$",
 i.e., data consisting of:
 \begin{itemize}
 \item An object $X_\rho\in\Ac$ for any $m$-dimensional face $\rho\subset\Delta^m$.
 \item An (iso)morphism 
 \[
 \phi_\sigma: \bigotimes_{i\equiv 0\mod 2} X_{\partial_i\sigma} \lra  \bigotimes_{i\equiv 1\mod 2}
 X_{\partial_i\sigma}
 \]
 for each $(m+1)$-dimensional face $\sigma\subset\Delta[n]$.
 \item For  any $(m+2)$-dimensional face $\tau\subset\Delta[n]$
 the isomorphisms $\phi_{\partial_i\tau}$ must satisfy a compatibility
 condition,  identical to one 
 spelled out in \cite{Previdi}, \S 2.9 (for the case when $\Ac$ is the category of
 torsors over an abelian group $A$). 
 \end{itemize}
 Note that for $m=0$ this gives the usual definition of the nerve of $\Ac$. Further details
 are left to the reader. 
 \qed
 
 \newtheorem*{examples-picard-spectra}{(3.1.7) Examples}
 \begin{examples-picard-spectra}
 (a) The Picard category $[A^\bullet]$ from Example (3.1.2)(c) corresponds to the
 Eilenberg-MacLane spectrum $EM(V^\bullet)$. 
 
 \vskip .2cm
 
 (b) Let  $Y$ be  any {\em connective}  spectrum , i.e., a spectrum such such that $\pi_i(Y)=0$ for $i<0$.
 Then we have a well defined truncation $Y_{[0,1]}$ of $Y$ in the homotopy degrees $0,1$,
 which is a $[0,1]$-spectrum and can therefore be described in terms of an appropriate Picard
 category. We consider several particular cases.
 
 \vskip .2cm
 
 (c) (K-theory spectra.) Let $\Ec$ be an exact category in the sense of Quillen
 and
 $\Kc(\Ec)$ be its algebraic K-theory spectrum, so $\pi_i\Kc(\Ec) = K_i(\Ec)$
 are the algebraic K-groups of $\Ec$. The symmetric Picard category corresponding to
 the $[0,1]$-spectrum $\Kc(\Ec)_{[0,1]}$ was described by Deligne in 
 \cite{Deligne:Determinant}. 
 It was called the {\em category of virtual objects} of $\Ec$
 and denoted $\Vc(\Ec)$. This Picard category can be characterized
 as the target of the universal determinantal theory on $\Ec$, see \cite{Previdi},
 Ex. 2.13, 2.29.

 \vskip .2cm
 
 (d) (K-theory of a field.) Let $\kk$ be a field and $\Ec$ be the category of finite-dimensional
 $\kk$-vector spaces. In this case
 $$K_0(\Ec) = K_0(\kk) = \ZZ, \quad K_1(\Ec) = K_1(\kk) = \kk^\times.$$
 It was observed in \cite{Deligne:Determinant} and empasized in 
   \cite{Drinfeld} that $\Vc(\Ec)$ is equivalent, as a symmetric Picard category,
   to the category $\Pc ic^\ZZ(\kk)$ from Example 3.1.2(d).

 \vskip .2cm
 
 (e) (The spherical spectrum.) Let $\SS$ be the spherical spectrum, so $\pi_i(\SS)=\pi_i^{\st}$
 are the stable homotopy groups of spheres:
 $$\pi_0^{\st}=\ZZ, \quad \pi_1^{\st} = \ZZ/2, \quad \pi_2^{\st} = \ZZ/2, \quad \pi_3^{\st} = \ZZ/24, 
 \text{  etc.} $$
 As $\SS$ is connective, we have the truncation $\SS_{[0,1]}$. The corresponding Picard
 category can be viewed as a free symmetric Picard category generated by one object. It can
 be described  as a modification of the category   $\Pc ic^\ZZ(\kk)$. 
  More precisely, let  $\Pc ic^\ZZ(\ZZ)$
 be the category of $\ZZ$-graded  abelian groups $L=\bigoplus_{n\in \ZZ} L_n$
 which are free of total rank 1. Morphisms are isomorphisms of graded abelian groups. 
 Tensor product and symmetry are defined similarly to $\Pc ic^\ZZ(\kk)$.

   In other words, the Picard category classifying (the truncation of)
   the spherical spectrum $\SS$
   can be seen as the ``sign skeleton"
   of the category of  graded vector spaces with Koszul sign rule: it contains exactly the
   minimal data necessary to write down this rule. So it is the structure of $\SS$
   which ultimately leads to the sign rules of super-mathematics. 
   
   \vskip .2cm
   
   (f) (The loop space of the spherical spectrum.) Remarkably, the spherical spectrum leads
   to super-constructions not in one, but in two ways. More precisely,  let
   $\SS_0$ be the connected component of $\SS$ corresponding to
   the zero element of $\pi_0(\SS)=\ZZ$. The loop space (spectrum) $\Omega \SS_0$ is then
   a connective spectrum whose homotopy groups are the $\pi_i^{\st}$ 
   but with shifted numeration and with
   $\pi_0^{\st}=\ZZ$ disregarded:
   $$
   \pi_0(\Omega\SS_0) = \ZZ/2, \quad \pi_1(\Omega\SS_0)=\ZZ/2, \quad \pi_2(\Omega\SS_0)=\ZZ/24, \text
   { etc.}
   $$ 
   Therefore the truncation $(\Omega\SS_0)_{[0,1]}$ gives rise to a symmetric
   Picard category $\Pc$
   with
   $$\pi_0(\Pc) = \pi_1(\Pc) = \ZZ/2.$$
  \end{examples-picard-spectra}
  
  \newtheorem*{thm:spherical-super-lines}{(3.1.8) Theorem}
  \begin{thm:spherical-super-lines}  As a symmetric
  Picard category,  $\Pc$ is equivalent to the category $\Pc ic^{\ZZ/2}(\ZZ)$, whose
  objects are  $\ZZ/2$-graded free abelian
  groups of total rank $1$, with the usual $\ZZ/2$-graded tensor product and
  the Koszul sign rule for the symmetry. 
    \end{thm:spherical-super-lines}
    
    In other words, the truncated spectrum $(\Omega\SS_0)_{[0,1]}$ can be identified with
    the ``reduction modulo 2" of $\SS_{[0,1]}$. 
    
    \vskip .2cm
    
    We will not prove Theorem (3.1.8) in this paper 
     but will use it as a topological motivation for our
    explicit constructions involving
    symmetric groups. Let us explain the connection in more detail. 
    
    \vskip .2cm
    
    \noindent {\bf (3.1.9) The higher sign character of $S_n$
    via the Barratt-Priddy-Quillen theorem.} The classical 
     Barratt-Priddy-Quillen (BPQ) theorem \cite{Priddy, Barratt:Priddy}
      gives a construction of the
     spherical spectrum $\SS$ in terms of the symmetric groups $S_n$,
     $n\geq 1$. More precisely, let $S_\infty = \varinjlim_n S_n$ be the
     union of the standard chain of embeddings and $BS_\infty$ be its
     classifying space. The alternating group $A_n\subset S_n$ is,
     for $n\geq 5$, equal to its commutant, and 
     the same is true for the union $A_\infty\subset S_\infty$.
     In this situation, one can form the {\em Quillen plus-construction}
      $BS_\infty\to B^+S_\infty$. It can be characterized uniquely
     (up to homotopy equivalence) as the map reducing $\pi_1(BS_\infty)=S_\infty$
     to $S_\infty/A_\infty=\ZZ/2$ and inducing an isomorphism on integral homology
     groups.  The BPQ theorem says that $B^+S_\infty=\Omega^\infty \SS_0$ 
     is the connected component of the infinite loop space corresponding to
     the spherical spectrum. Further, the spectrum structure (i.e., the deloopings
     of $B^+S_\infty$) can be canonically recovered using the semigroup
     structure on $B^+S_\infty$ coming from the direct sum embeddings
     $S_m\times S_n\to S_{m+n}$.
     
     The  BPQ theorem implies that we have a canonical map
     \[
     BS_n\lra BS_\infty \lra B^+S_\infty \sim \Omega^\infty \SS_0.
     \]
     Passing to the loop spaces, we get a map
     \[
   \mathfrak {sgn}: S_n\lra \Omega^\infty(\Omega\SS_0)
     \]
     which we want to call the {\em homotopy-theoretic sign character}. 
     This map and its (higher) categorical manifestations should be the
     right tool for studying  analogs of exterior powers in
     the context not only of categories, as in this paper, but of
     higher categories as well. 
     
     In particular, replacing  the spectrum $\Omega\SS_0$ by its $[0,1]$-truncation
     and then by the Picard category $\Pc$ as in Example (3.1.7)(f), we get a canonically
     defined Picard character
     \[
     \operatorname{Sgn}: S_n \lra \Pc,
     \]
 which we call the {\em categorical sign character} of $S_n$. On the level of $\pi_0$,
 it induces the ordinary sign character $\sgn$. 
 
 Combining the above with Theorem (3.1.8),  we get  a {\em canonically defined}
   Picard character
 of $S_n$ with values in the Picard category $\Pc ic^{\ZZ/2}(\ZZ)$ of
 super-lines over $\ZZ$. As $\Pc ic^{\ZZ/2}(\ZZ)$ naturally acts on any super-linear
 category,  we get a conceptual explanation of the
 following well known phenomenon: 
  the theory of projective representations of symmetric
 groups  becomes much streamlined and simplified by systematic
 use of super-objects. One can find many instances of this
 phenomenon in \cite{Kleshchev}. 
 In the following sections we define 
   the Picard character $\operatorname{Sgn}: S_n\to\Pc ic^{\ZZ/2}(\ZZ) $ 
   in an elementary way and apply it to the construction of exterior
   powers of categories.

\vskip .3cm

\subsection*{(3.2) The central extension of $S_n$. Naive exterior powers.} 
Assume $\on{char}(\kk)\neq 2$. 
The sign character 
(1-dimensional representation) $\sgn: S_n\to\{\pm
1\}$ generates the group
$$H^1(S_n, \kk^*) \quad = \quad \on{Hom}(S_n, \kk^*) \quad = \quad \ZZ/2.$$
 So the most straightforward way of generalizing it to the categorical case would be to
 use 1-dimensional
2-representations of $S_n$, which correspond (Examples 2.1.2) to elements of
$H^2(S_n, \kk^*)$. It is well known, see \cite{Kleshchev},
that
$$H^2(S_n, \kk^*) = \ZZ/2, \quad n\geq 4. $$
This cohomology group is generated by the class of the central extension
$$1\to \{\pm 1\} \to \widetilde{S}_n \buildrel\varpi\over\longrightarrow  S_n \to 1,
\leqno (3.2.1)$$
which makes sense for any $n\geq 1$ and was first found by I. Schur
\cite[p.164]{Schur}. To describe $\widetilde S_n$, 
 recall that $S_n$ is generated by elementary
transpositions
$$\sigma_i = (i, i+1), \quad i=1, ..., n-1,\leqno (3.2.2)$$
subject to the relations
$$\sigma_i^2=1, \quad \sigma_i\sigma_{i+1}\sigma_i = \sigma_{i+1}\sigma_i\sigma_{i+1},
\quad \sigma_i\sigma_j = \sigma_j \sigma_i, \quad |i-j|\geq 2. 
\leqno (3.2.3)$$ 
Now, $\widetilde{S}_n$ is generated by elements $s_1, ..., s_{n-1}, \zeta$
subject to the following relations:
$$\zeta^2=1, \quad \zeta s_i = s_i\zeta; \leqno (3.2.4)(a)$$
$$s_i^2=\zeta, \quad 
s_is_{i+1}s_i = s_{i+1} s_i s_{i+1}, \quad s_is_j =\zeta s_j s_i, \quad |i-j|\geq 2.
\leqno (3.2.4)(b)$$
The central subgroup $\{\pm 1\}$ in $\widetilde{S}_n$ is generated by $\zeta$. 

\vskip .2cm

The extension $\widetilde{S}_n$ of $S_n$ by $\{\pm 1\}$ gives, after the change of groups
$\{\pm 1\} \hookrightarrow \kk^*$, an extension by $\kk^*$. Let us describe
this extension in the spirit of Example 2.1.2(b), i.e., in terms of 1-dimensional
$\kk$-vector spaces
$$L_\sigma \quad = \quad \bigl( \varpi^{-1}(\sigma) \otimes_{\{\pm 1\}} \kk^*\bigr) \cup \{0\}, 
\quad\sigma\in S_n. \leqno (3.2.5)$$
Let $D(\sigma)$ be the set of all reduced decompositions (i.e., decomposition of the
minimal possible length $l(\sigma)$)
$$\sigma = \sigma_{i_1} ... \sigma_{i_l}, 
\quad l=l(\sigma). \leqno (3.2.6)$$
We view an element of $D(\sigma)$ as a sequence of indices
$$\bold{i} = (i_1, ..., i_l), \quad i_\nu \in\{ 1, ..., n-1\}. $$
Any two decompositions $\bold{i}, \bold{i}'
\in D(\sigma)$ can be transformed into each other by a sequence of moves
of two types: {\em hexagonal moves}
$$\bold{i} = (..., i, i+1, i, ...) \quad \rightsquigarrow\quad (..., i+1, i, i+1, ...) = \bold{i}',$$
and {\em square moves}
$$\bold{i} = (..., i, j, ...) \quad\rightsquigarrow\quad (..., j, i, ...) = \bold{i}', 
\quad |i-j|\geq 2.$$
These moves correspond to implementing the second and third type of relations in (3.2.3).

\newtheorem*{237}{(3.2.7) Proposition}
\begin{237}
(a) For any sequence of moves transforming $\bold{i}$ into $\bold{i}'$ the number of square
moves has the same parity.

(b) The line $L_\sigma$ is identified with the space of functions $f: D(\sigma)\to\kk$ which are
unchanged under hexagonal moves and change sign under each square move.
\end{237}

\noindent {\sl Proof:} Part (b) is a direct translation of the definition
of $\widetilde{S}_n$ by generators and relations (3.2.4). Indeed, any $\bold{i}\in D(\sigma)$
gives a lift
$$\widetilde{\sigma}_{\bold{i}} = s_{i_1} ... s_{i_l} \in\widetilde{S}_n, \quad 
\varpi(\widetilde{\sigma}_{\bold{i}})=\sigma, \leqno (3.2.8)$$
and these lifts behave as claimed under the two types of moves. 

Part (a) is an implicit preliminary to (b), expressing the fact that 
nonzero $f$ as in (b) exist. It is equivalent to the (classical)
fact that the group defined by the relations (3.2.4) does not collapse,
i.e., is indeed an extension as in (3.2.1).  \qed.

\vskip .2cm

Let us choose a reduced decomposition for each $\sigma\in S_n$. Then we get
a lifting $\widetilde{\sigma}$ for each $\sigma$, so we have a 2-cocycle
$$c: S_n\times S_n\to \{\pm 1\}, \quad c(\sigma, \tau) = \widetilde{(\sigma\tau)}
(\widetilde{\tau})^{-1} (\widetilde{\sigma})^{-1}, \leqno (3.2.9)$$
describing the extension.

\vskip .2cm

\newtheorem*{2.3.10}{(3.2.10) Definition}
\begin{2.3.10}
(a) The naive sign 2-representation $\on{Sgn}_{\on{naive}}$ of $S_n$
is the 1-dimensional 2-representation $\mV_c$, see Example 2.1.2(a),
corresponding to $c$ from (3.2.9).

(b) Equivalently, $\on{Sgn}_{\on{naive}}$ can be described as the action of
$S_n$ on $\on{Vect}^f$ by the functors
$$\varrho(\sigma): V\mapsto L_\sigma\otimes V, \quad \sigma\in S_n,$$
and $\phi_{\sigma,\tau}$ being induced by the multiplicativity
isomorphisms $L_\sigma\otimes L_\tau\to L_{\sigma\tau}$. 
\end{2.3.10}

%The reason for the term ``naive'' will be explained later.
The following is a categorical analog of the main construction in
\cite{Dijkgraaf}. 

\newtheorem*{2.3.11}{(3.2.11) Definition}
\begin{2.3.11}
Let $\mV$ be a linear category. The naive exterior power of $\mV$ is defined to be
the category
$$\Lambda^n_{\on{naive}}(\mV) \quad =\quad \bigl( \mV^{\boxtimes n} \boxtimes
\on{Sgn}_{\on{naive}}\bigr)^{S_n}.$$
\end{2.3.11}
By restricting to the generators of $S_n$, we get the following:

\newtheorem*{2.3.12}{(3.2.12) Reformulation}
\begin{2.3.12}
Explicitly, an object of $\Lambda^n_{\on{naive}}(\mV)$ consists of an object
$V\in\mV^{\boxtimes n}$ together with isomorphisms
$$s_i: V\to \sigma_i^*(V), \quad i=1, ..., n-1,$$
satisfying the following conditions:
\vskip .1cm

(a) For any $i = 1, ..., n-1$, the composition
$$V\buildrel s_i\over\longrightarrow \sigma_i^*(V) 
\buildrel\sigma_i^*(s_i)\over\longrightarrow \sigma_i^*\sigma_i^*V 
\buildrel\simeq \over\longrightarrow V$$
is equal to $(-\id_V)$. 

(b) For any $i=1, ..., n-2$ the two compositions
$$V\buildrel s_i\over\longrightarrow \sigma_i^*V
\buildrel\sigma_i^*(s_{i+1})\over\longrightarrow \sigma_i^*\sigma_{i+1}^* V
\buildrel \sigma_i^*\sigma_{i+1}^*(s_i) \over\longrightarrow \sigma_i^*\sigma_{i+1}^*\sigma_i^* V,$$

$$V\buildrel s_{i+1}\over\longrightarrow \sigma_{i+1}^*V
\buildrel\sigma_{i+1}^*(s_{i})\over\longrightarrow \sigma_{i+1}^*\sigma_{i}^* V
\buildrel \sigma_{i+1}^*\sigma_{i}^*(s_{i+1}) \over\longrightarrow \sigma_{i+1}^*\sigma_{i}^*\sigma_{i+1}^* V$$
become equal after identifying the rightmost terms with
$$(\sigma_i\sigma_{i+1}\sigma_i)^*V = (\sigma_{i+1}\sigma_i\sigma_{i+1})^*V.$$

(c) For any $i, j$ such that $|i-j|\geq 2$ the two compositions
$$V\buildrel s_i\over\longrightarrow \sigma_i^*V \buildrel \sigma_i^*(s_j)\over\longrightarrow
\sigma_i^*\sigma_j^*V, \quad \quad
V\buildrel s_j\over\longrightarrow \sigma_j^*V \buildrel \sigma_j^*(s_i)\over\longrightarrow
\sigma_j^*\sigma_i^*V $$
differ by sign, after identifying the rightmost terms with
$$(\sigma_i\sigma_j)^*V = (\sigma_j\sigma_i)^*V.$$

\end{2.3.12}

\vskip .3cm

\noindent {\bf (3.3) Projective representations of $S_n$. Role of superalgebra.}
Let $c$ be a cocycle defining $\widetilde{S}_n$, as in (3.2.9). Denote by
$\on{Rep}_c(S_n)$ the category of finite-dimensional projective representations of
$S_n$ with central charge $c$, i.e., of vector spaces $V$ together with a map
$\varrho: S_n\to \on{Aut}(V)$ satisfying
$$\varrho(\sigma\tau) = c(\sigma, \tau) \cdot \varrho(\sigma) \varrho(\tau).$$
More invariantly, an object of $\on{Rep}_c(S_n)$ is a representation of $\widetilde{S}_n$
with $\zeta$ acting by $(-1)$. By passing to the generators of $\widetilde{S}_n$,
we see that an object of $\on{Rep}_c(S_n)$ is the same as $V\in\on{Vect}^f$ together
with automorphisms $s_i: V\to V$ satisfying the relations
$$s_i^2=-1, \quad s_is_{i+1}s_i = s_{i+1}s_i s_{i+1}, \quad s_is_j = - s_j s_i, \quad |i-j|\geq 2.
\leqno (3.3.1)$$
Let $\kk[S_n]^c$ be the associative $\kk$-algebra with generators $s_i$ subject
to these relations. It is the quotient of the group algebra of $\widetilde S_n$ by the
relation $\zeta = -1$. Then we can say that $\on{Rep}_c(S_n)$ is identified with
the category of finite-dimensional $\kk[S_n]^c$-modules. 

\vskip .2cm

\noindent {\bf (3.3.2) Example.} Taking $\mV = \on{Vect}^f = [1]$, as in Example 2.2.2, we find, by
comparing (3.2.12) with (3.3.1), that
$$\Lambda^n_{\on{naive}}(\mV) \quad = \quad \on{Rep}_c(S_n). $$

\vskip .2cm

Irreducible objects of $\on{Rep}_c(S_n)$ in characterisic 0 were described by Schur \cite{Schur}.
We recall his results. 

\vskip .2cm

\newtheorem*{243}{(3.3.3) Definition}
\begin{243}
A {\em strict partition} of $n$ is a partition
$$n=\sum_{k\geq 1}kN_k,$$ such that all the coefficients $N_k$ are either
$0$ or $1$. The number of strict partitions of $n$ will be denoted
$s(n)$. 
A partition of $n$
is called {\em even} (respectively {\em odd}) if the corresponding
conjugacy class in $S_n$ is even (respectively odd). 
\end{243}

\vskip .2cm

\newtheorem*{244}{(3.3.4) Theorem (Schur)}
\begin{244} Assume $\on{char}(\kk)=0$. Then:

(a) Irreducible objects of $\on{Rep}_c(S_n)$ are labelled by data
consisting of a strict partition $\lambda$ of $n$ together with, if $\lambda$ is
odd, a sign $+$ or $-$. Thus for an even $\lambda$ there is one irreducible
representation $V_\lambda$, while an odd $\lambda$ gives rise to
two representations $V^+_\lambda$, $V^-_\lambda$. 

(b) If $\lambda$ is even then 
$$V_\lambda\otimes\sgn = V_\lambda,$$
while for $\lambda$ odd we have
$$V^+_\lambda \otimes\sgn = V^-_\lambda.$$\qed
\end{244}

Note that the tensor product of an object of $\on{Rep}_c(S_n)$ and an
actual representation of $S_n$ is again an object of
$\on{Rep}_c(S_n)$, thus giving sense to (b). 

\vskip .2cm

\noindent {\bf (3.3.5) Remarks.} (a) If we understand the term
``projective representation'' in the classical sense, as a homomorphism
to $PGL_N$ for some $N$, then $V_\lambda^+$ and $V_\lambda^-$
give the same homomorphism. So, conjugation classes of irreducible
homomorphisms $S_n\to PGL_N$ not lifting to $GL_N$ are in bijection
with strict partitions of $n$. This is the original formulation of
Schur (\cite{Schur}, p. 156). 

\vskip .2cm

(b) Note also that the   generating function for strict partitions
$$\sum_{n\geq 0} s(n) q^n \quad = \quad \prod_{n\geq 1} (1+  q^n) \leqno (3.3.6)$$
reminds of $\phi(q)$, which is the inverse of the generating function
for 
$$p(n) \quad = \quad \on{Dim} \,\, \Sym^n (\on{Vect}^f). $$
This is  the type of relation existing between the generating 
functions for exterior and symmetric powers in usual linear algebra:
$$ \( \sum \, \dim \,\Lambda^n(V) q^n\) \cdot \( \sum  \, \dim S^n(V) (-q)^n\) 
\quad = \quad 1.\leqno (3.3.7)$$
This suggests that one should modify the ``naive'' definitions
of $\on{Sgn}$ and $\Lambda^n(\mV)$ above so as to get
$\Lambda^n(\on{Vect}^f)$   satisfying an analog of (3.3.7). 

\vskip .2cm

It turns out that
the key to such a modification is provided by the superalgebra
point of view \cite{Jozefiak:Projective} on projective representations of $S_n$. 
 It begins with the observation that the third group of relations
in (3.3.1) (anti-commutativity of the $s_i$) can be
seen  as an instance of super-commutativity.
For this to hold, we have to consider the twisted group algebra $\kk[S_n]^c$
as a super-algebra, i.e. equip it with the $\ZZ/2$-grading defined by
$\deg(s_i) = \bar 1\in\ZZ/2$. We will now develop  
some additional features of this approach.

\vskip .3cm

\noindent {\bf (3.4) Super 2-representations}

\newtheorem*{251}{(3.4.1) Definition}
\begin{251}
  By a super 2-representation of a group $G$, we will mean an action
  of $G$ on a superlinear category $\mV$ by superlinear functors and
  supernatural transformations.
\end{251}
\newtheorem*{252}{(3.4.2) Example}
\begin{252}
By (1.6.3), a superlinear action of $G$ on $\SVect$ 
is the same as a Picard character $\Lambda: G\to\Pc ic^{\ZZ/2}(\kk)$,
i.e., a family of super-lines $\Lambda_g, g\in G$, equipped with
compatible multiplicativity  and unit isomorphsms
$$\nu_{g,h}: \Lambda_g\otimes\Lambda_h\lra \Lambda_{gh}, \quad
\nu_1: \Lambda_1\lra\kk,$$
  see Definition 3.1.3.
 In particular, $G$ acquires a $\ZZ/2$-grading
 \[
  G\longrightarrow \ZZ/2\quad 
  g\longmapsto \bar g:= \on{deg}(\Lambda_g).
 \]
Up to non-canonical equivalence of 2-representations, we may assume
that $\Lambda_g = \Pi^{\bar g}\kk$, that $\nu_1=\id$ and that the
$\nu_{g,h}$ are multiplication by the values $c(g,h)\in\kk^*$ of a
normalized 2-cocycle of $G$. 
We will refer to a 2-representation in this form as
the {\em cocycle super
  2-representation}  $\varrho_c$. 
  
  Note the particular case $G=S_n$. In this case  the sign 
   homomorphism $\sgn: S_n\to \{\pm 1\}$ can be seen as a $\ZZ/2$-grading
  \[
  S_n\lra \ZZ/2, \quad \sigma\longmapsto \bar\sigma, \quad (-1)^{\bar\sigma}=\sgn(\sigma).
  \]
\end{252}
\newtheorem*{253}{(3.4.3) Definition}
\begin{253}
  The sign super 2-representation $\varrho_{\on{sgn}}$ of $S_n$ is the cocycle
  super 2-representation corresponding to the
  central extension $\widetilde S_n$ of $S_n$ and to the $\ZZ/2$-grading corresponding
  to the sign character,
  i.e.,
  $$
    \Lambda_\sigma := \Pi^{\bar\sigma} L_\sigma,
  $$
  and
  $\nu_{g,h}$ is defined by multiplication in $\widetilde S_n$ (see
 (3.2.5) ff.). 
\end{253}

\vskip .3cm
\noindent {\bf(3.4.4) Characters of cocycle super 2-representations.}
We now find the categorical character and the 2-character of
a cocycle 2-representations $\varrho_c$ as in Example (3.4.2). 
For a 2-representation
$\varrho$ of a group $G$ on a linear category $\mV$, its
categorical character is the family of $\kk$-vector
spaces
$$X_\varrho(g)\medspace  = \medspace 
\on{\mathcal N\!at}(\id_\mV, \varrho(g)), \quad g\in G,
$$
see \cite{Ganter:Kapranov}. These vector spaces form a conjugation
invariant vector bundle on the discrete set $G$, i.e., we have
isomorphisms
$$\isomap{\psi_s}{X_\varrho(g)}{X_\varrho(sgs^{-1})},
$$
compatible with compositions of the $h$'s, see \cite[Prop
4.10]{Ganter:Kapranov}.  
These constructions are extended to the superlinear case as follows:
\newtheorem*{255}{(3.4.5) Definition}
\begin{255}
Let $\varrho$ be a super-linear 2-representation. Then 
the categorical (super-) character of $\varrho$ is the super-vector space
$$sX_\varrho^{\bullet}(g) \medspace =\medspace 
\on{s\mathcal N\!at_\bullet}(\id_\mV, \varrho(g)).$$  
To define the action of 
$\psi_s$ on $sX^{\bar1}$, we note that superlinearity of $\rho$ makes 
$\Pi$ an equivariant functor whose flip isomorphisms
$$
  \tau_s\negmedspace:\rho(s)\Pi\stackrel\cong\longrightarrow
  \Pi\rho(s)
$$
are given by the superlinearity data of the $\rho(g)$.
Supressing $\varrho$ from the notation, $\psi_s$ sends
$\xi\in sX^{\bar1}(g)$ to the composite\footnote{This should be compared
  to the string diagram at the beginning of Section 5.3 in 
\cite{tracesofsymmetricpowers}.}
$$
  (\Pi\phi_{s,g,s\inv})\circ(\tau_s gs\inv)\circ(s\xi
  s\inv)\circ\phi_{s,s\inv}\inv\circ\phi_1\inv. 
$$
Assuming all the $sX_\varrho^\bullet(g)$ are finite-dimensional, the
2-character of a super-linear 2-representation $\varrho$ is 
defined as
$$\chi_\varrho(g,h) \quad = \quad \on{tr} \bigl\{ \psi_g(h): sX_\varrho^\bullet(g)
\to sX_\varrho^\bullet(g)\bigr\}, \quad gh=hg.$$ 
This is a $\kk^2$-valued function, defined on pairs of commuting
elements of $G$.  
\end{255}
Let now $c$ be a cocycle on a $\ZZ/2$-graded group $G$, and let 
$$
  1\to\kk^*\to\widetilde G\to G\to 1
$$
be the central extension classified by $c$. For a pair of commuting
elements $g,h\in G$, we let 
$$
  \epsilon(g,h) := \frac{c(g,h)}{c(h,g)} 
$$
and
$$
  e(g,h) :=
  (-1)^{\on{deg}(g)\on{deg}(h)}\frac{c(g,h)}{c(h,g)}. 
$$
In terms of the group extension, $\epsilon(g,h)$ is the {\em symbol}
of $h$ and $g$ in the sense of cf \cite[\S 8]{Milnor}. More precisely, for any
lifts $\widetilde g$ of $g$ and $\widetilde h$ of $h$, we have
$$
  \epsilon = 
  \widetilde h\inv\widetilde g\inv\widetilde h\widetilde g.
$$ 
We will refer to $e(g,h)$ as the {\em supersymbol} of $g$ and $h$.
\newtheorem*{256}{(3.4.6) Proposition}
\begin{256}
The 2-character of $\varrho_c$ is given by
$$
  \chi_{\varrho_c}(g,h) = 
  \begin{cases}
    (e(g,h)\mid 0) & \text{if $g$ is even and}\\
    (0\mid e(g,h)) & \text{if $g$ is odd.}
  \end{cases}
$$
\end{256}

\begin{proof}
By definition of $\varrho_c$, its
categorical character is given by
$$sX_{\varrho_c}^\bullet(g) = \Lambda_g = \Pi^{\bar g} L_g.$$
So,  
$$
\chi_{\varrho_c}(g, g) = 
  \begin{cases}
    (\psi_h\mid0) & \text{if g is even and}\\
    (0\mid\psi_h) & \text{if g is odd.}\\
  \end{cases}
$$
For even $g$, the map $\psi_h$ is the 
scalar by which the conjugation by $h$
acts on the 1-dimensional vector space $L_g$. From the definition
(3.2.5) of $L_g$ we see that $\psi_g(h)$ describes
the way by which the conjugation by (any lift of) $h$ in $\widetilde
G$ acts on the set $\varpi^{-1}(g)$.
Hence $\psi_g(h) = \epsilon(g,h)$, as claimed. 
For odd $g$, the definition of $\psi_h$ in Definition (3.4.5) includes
the flip map $\tau_h$. By the Koszul sign rule in (1.6.3), $\tau_h$ is
multiplication with $(-1)^{\bar h}$. The rest is as in the even case
(compare \cite[Lem.6.1]{innerproducts}).
\end{proof}

\vskip .3cm

\noindent {\bf (3.4.7) The character of the sign super 2-representation.}
In the case of $\widetilde S_n$, the symbol $\epsilon(\sigma,\tau)$
was calculated by Dijkgraaf in \cite{Dijkgraaf}.
%Note, however, that there is a dicrepancy between Schur's
%results\footnote{The relevant theorem is \cite[Thm IV]{Schur},
%  see (2.6.9) below.}
%and the main result of \cite{Dijkgraaf}: 
%\marginpar{Question 10}
%according to Schur,
%I do not know whether this makes any sense from a
%physics point of view, but indeed, the first equality of Formula (3.20) in
Note that in Dijkgraaf's calculation, Formulas (3.20)ff. 
should include the sign
%is wrong: the analogous step for (3.19) uses
%that $x_2$ has even parity. This argument breaks down for odd $n$, and
%instead Formula (3.20) needs to include the sign
$$
  \epsilon(x_n,g') = (-1)^{n|g'|}.
$$ 
%In other words, whether the odd long strings 
%behave bosonically or fermionically depends on whether the number of
%even long strings is even or odd. Whenever you add an even long
%string, the odd long strings change their behaviour from bosonic to
%fermionic or back. 
%
%\newtheorem*{2512}{(2.5.12) Example}
%\begin{2512}
%Let
%\begin{eqnarray*}
%\sigma&=&(123)(456)(78)\\
%k'&=&(123)\\
%5k''&=&(456)\\
%x_3&=&(14)(25)(36)\\
%\tau:= x_3k''& =& (142536).  
%\end{eqnarray*}
%Then $|\sigma|=1=|\tau|$, and $\sigma=\tau^2(78)$. Dijkgraaf calculates
%$$
%  \epsilon(\sigma,\tau) = \epsilon(\sigma,x_3)\epsilon(\sigma,k'') = 1,
%$$ 
%by \cite[(3.14), (3.20)]{Dijkgraaf}. On the other hand, we have
%$$
%  \epsilon(\sigma,\tau) = \epsilon(\tau^2,\tau)\cdot \epsilon((78), \tau) =
%  1\cdot(-1). 
%$$  
%\end{2512}
With this correction, Dijkgraaf's result allows us to calculate
the 2-character of $\varrho_{sgn}$. 
\newtheorem*{258}{(3.4.8) Definition}
\begin{258}
  Let $H$ be a finite group, and let $\map\varrho HGL(V_\bullet)$ be a
  superrepresentation of 
  $H$. Then we have the representation
  $$
    \longmap{\varrho\wr\sgn}{H\wr\Sn}{GL(V^{\tensor n}_\bullet)},
  $$
  where $H^n$ acts by $\varrho^{\tensor n}$, and $\Sn$ permutes the tensor
  factors, adding a sign.
  Further, for even $k$, we let $\varrho_{\on{spin}}$ be the
  $(0|1)$-dimensional superrepresentation of
  $\ZZ/k\ZZ$ with character $(0|\chi_{\on{spin}})$, where
  $$
    \longmap{\chi_{\on{spin}}}{\ZZ/k\ZZ}{\{\pm 1\}}
  $$
  is the unique non-trivial map.
\end{258}
%\marginpar{Question 11}
\newtheorem*{259}{(3.4.9) Theorem}
\begin{259}
  Let $\sigma\in\Sn$ be classified by the partition
  $n=\sum_{k\geq 1}kN_k$.
  Then the supersymbol
  $e(\sigma,\tau)$, viewed as (one-dimensional)
  superrepresentation 
  of the centralizer $$C_\sigma\cong\prod_{k\geq1}\ZZ/k\ZZ\wr
  S_{N_k},$$ 
  equals
  $$
    e(\sigma,\tau)= \bigotimes_{k\text{
        even}}\chi_{\on{spin}}\wr\on{sgn}_{S_{N_k}}.$$
%  if $\sigma$ is even and
%  $$  
%    \epsilon(\sigma,\tau)=\bigotimes_{k}\on{sgn}_{S_{N_k}}
%  $$
%  if $\sigma$ is odd.
\end{259}
Together with Proposition (3.4.6) this determines the 2-character 
$\chi_{\on{sgn}}$ of $\varrho_{\on{sgn}}$.  

\vskip .3cm

%\vfill\eject 

\noindent {\bf (3.5) Equivariant objects.} 
Throughout this section, we assume the characteristic of $\kk$ to be zero.
Let
$\varrho$ be a super-2-representation of $G$ on $\mV$.
Consider the category $\mV^G$ of equivariant objects in $\mV$.
If $\Pi$ is $G$-equivariant,%\marginpar{is that always the case?} 
then the category $\mV^G$ is itself a
superlinear category, with 
$$
  \Pi(V,\{\epsilon_g\}) := \(\Pi(V),((\tau_g)_V)\circ(\Pi\epsilon_g)\).
$$
Here ${\tau_g}\negmedspace :{\Pi g}\Rightarrow{g\Pi}$ is the natural
isomorphism making $\rho(g)$ a superlinear functor.
The graded extension $\(\mV^G\)_\bullet$ of $\mV^G$  
is canonically identified with the
full subcategory of $\(\mV_\bullet\)^G$ whose objects have 
only even structure maps.

\newtheorem*{261}{(3.5.1) Example}
\begin{261}
  Let $\mV=\SVect$, and let $G$ be a $\ZZ/2$-graded group acting on
  $\mV$ via $\varrho_c$ as in Example (2.5.2). Corresponding to the
  cocycle $c$, we have the 
  twisted group superalgebra $\kk[G^{\on{op}}]^c$. This is the
  associative $\kk$-algebra generated by
  $\{e_g\mid g\in G\}$ with multiplication
   $$e_g\cdot e_h = c(g,h)e_{hg}.\leqno (3.5.2)$$  
  It is made into a superalgebra, i.e., equipped with the $\ZZ/2$-grading,
   by putting $\deg(e_g)=\bar g$. 
  The category of equivariant objects, $\mV^G$, is equivalent, as a
  superlinear category, to the category $\kk[G^{\on{op}}]^c\on{-Smod^f}$ 
  of finite-dimensional left $\kk[G^{\on{op}}]^c$-supermodules and even
  morphisms between them. Note that the action of an odd element $g$
  on $\Pi V$ acquires the sign $\tau_g=-1$, as it should for left
  supermodules. 
\end{261}
%\newtheorem*{262}{(2.6.2) Lemma}
%\begin{262}
%  Assume that $\kk$ has characteristic prime to $|G|$. Then the
%  forgetful map
%  $$
%    \longmap{f}{\on{Center}\(\(\mV_\bullet\)^G\)}
%    {\on{sCenter_\bullet}\(\(\mV^G\)\)}
%  $$
%  is an isomorphism.
%\end{262}
%\begin{proof}
%  Let $(V,\{\epsilon_g\})$ be an object of $(\mV_\bullet)^G$, where we do
%  not assume that the $\epsilon_g$ are even maps. Let 
%  $$
%    AV = \bigoplus_{g\in G}\varrho(g)V.
%  $$
%  The isomorphisms
%  $\map{\phi_{s,g}}{\varrho(s)\varrho(g)}{\varrho(sg)}$ make 
%  $AV$ into an even equivariant object. We have the maps 
%  $$
%    d:=\longmap{\(\epsilon_g\)_{g\in G}}{V}{AV}
%  $$
%  and 
%  $$
%    a:=\longmap{\sum_{g\in G}\epsilon_g\inv}{AV}V
%  $$
%  in $\(\mV_\bullet\)^G$.
%  Their compositite equals $|G|\cdot\id_{(V,\{\epsilon_g\})}$. It
%  follows that $f$ is injective: let $\zeta$ be in the center of
%  $\(\mV_\bullet\)^G$. We will abbreviate
%  $\zeta_{(V,\{\epsilon_g\})}$ by $\zeta_V$ and 
%  $\zeta_{(AV,\phi_{g})}$ by $\zeta_{AV}$. Then $\zeta_V$ is
%  determined by $\zeta_{AV}$ via
%  $$
%    \zeta_V = \GG G a\circ\zeta_{AV}\circ d.     
%  $$
%  On the other hand,
%  given $\vartheta$ in the center of $\(\mV^G\)_\bullet$, this last formula
%  yields a well-defined extension of $\vartheta$
%  to all of $\(\mV_\bullet\)^G$. Hence $f$ is surjective.
%\end{proof}
Applying\footnote{The theorem is applied inside the 2-category
of super-linear categories, super-linear functors
and supernatural transformations.} 
\cite[Thm.5.13]{innerproducts},
we obtain
isomorphisms
$$
  \on{sCenter}\(\mV^G\) \cong\bigoplus_{g\in
    G}sX_\varrho(g)^{C_g}.
  \leqno (3.5.3)
$$
and
$$
  s\ttr(\Pi,\id)
  \cong\bigoplus_{g\in G}s\ttr\((\Pi,\id)\varrho(g)\)^{C_g}.
  \leqno (3.5.4)
$$
For the remainder of this Section, we specialize to the situation of
Example (2.5.2). So, 
we have a $\ZZ/2$-graded group
$G$, acting on $\SVect$ via a cocycle $c$.

Then the category $(\SVect)^G$ is semisimple with finitely many isomorphism
classes of simple objects.
%\marginpar{conditions? reference? Sergeev?} 

\newtheorem*{265}{(3.5.5) Definition}
\begin{265}
  A conjugacy class $[g]\sub G$ is called {\em $c$-regular} if 
  $$
    \(\forall h\in C_g\)\quad \(c(g,h) = c(h,g)\).
  $$ 
\end{265}
The following Corollary of (3.5.3) and (3.5.4) is a generalization of
\cite[Lem 2.11]{Brundan:Kleshchev}. 
\newtheorem*{266}{(3.5.6) Corollary}
\begin{266}
  The number of isomorphism classes of irreducible objects in
  $\on{Smod-}\kk[G]^c$, up to shift of grading, equals the number of
  even $c$-regular conjugacy classes in $G$. Among these isomorphism
  classes, the number   
  of self-associate ones equals the number of odd
  $c$-regular conjugacy classes in $G$. 
\end{266}
\begin{Pf}{}
  By (1.6.1)(c), the dimension of (3.5.3) counts the number of
  isomorphism classes of irreducible objects in $\mV^G$, up to shift of
  grading, while the dimension of (3.5.4) counts the number of self-associate irreducible
  objects. 
  For even $g$, the centralizer 
  $C_g$ acts by $\epsilon(g,-)$ on the one dimensional 
  vectorspace $sX_{\varrho_c}(g)$, and the $g$th summand of
  (2.6.4) is $s\ttr(\Pi,\id))=0$. Let $g$ be odd. Then
  $sX_{\varrho_c}(g)=0$, while 
  the $g$th summand of (3.5.4) equals $s\ttr(\id,-\id)\cong \kk$. 
  By (1.6.1)(c), the action of
  $h\in C_g$ differs by the sign $(-1)^{\bar h}$ from that in (3.5.6).
  In other words, $h$ acts again by multiplication with
  $\epsilon(g,h)$.     
\end{Pf}{}

We now specialize to the case of the symmetric group $\Sn$.
\newtheorem*{267}{(3.5.7) Theorem (Schur) {\cite[IV]{Schur}}}
\begin{267}%[{Schur \cite[IV]{Schur}}]
  An even conjugacy class $[\sigma]$ of $S_n$ is $c$-regular if and
  only if all cycles of $\sigma$ have odd length. An odd
  $[\sigma]$ is $c$-regular if and only if $\sigma$ possesses no two
  cycles of equal length $\geq 1$.  
\end{267}
\newtheorem*{268}{(3.5.8) Corollary}
\begin{268}[compare {\cite[Thm. 22.2.1]{Kleshchev}}]
  The number of isomorphism classes of irreducible objects of
  $\on{Smod-}\kk[\Sn]^c$, up to shift of grading, equals the number of
  strict partitions of 
  $n$. Among these, the number of absolutely irreducibles is equal
  to the number of even strict partitions, while the number of
  self-associates is equal to the number of odd strict partitions.
\end{268}
\begin{proof}{}
  By Euler's theorem, the number of strict partitions of $n$ equals the
  number of partitions of $n$ into odd summands. 
  The Corollary follows
  immediately.  
\end{proof}{}

\vskip .3cm

\noindent {\bf (3.6) Exterior powers of categories.} We now
modify the naive definition of  the exterior powers
by incorporating the super point of view.

\vskip .1cm

\newtheorem*{271}{(3.6.1) Definition}
\begin{271}
Let $\mV$ be a superlinear category. 
The $n$th exterior power of $\mV$ is defined as 
$$\Lambda^n(\mV) \quad = \quad \bigl( \mV^{\boxtimes_s n} \boxtimes_s \on{Sgn}
\bigr)^{S_n}.$$
If $\mV$ is an abelian superlinear category
and $\mV^{\widehat\boxtimes_s n}$ exists, 
or if $\mV$ is a 2-periodic perfect dg-category, we define  
$$\widehat{\Lambda}^n(\mV) \quad = \quad \bigl(
\mV^{\widehat{\boxtimes}_s n} \boxtimes_s \on{Sgn} 
\bigr)^{S_n}.$$
There is no need to complete the super tensorproduct with
$\on{Sgn}$. 

\vskip .2cm
Explicitly, the exterior power
$\Lambda^n(\mV)$ has, as objects, data consisting of: 

\vskip .1cm

\noindent (a) An object $V\in \mV^{\boxtimes_s n}$; 

\vskip .1cm  

\noindent (b) Isomorphisms
$$s_i: V\to \Pi \sigma_i^*(V), \quad i=1, ..., n-1,$$
satisfying the following conditions:

\vskip .1cm

\noindent (c)  For any $i = 1, ..., n-1$, the composition
$$V\buildrel s_i\over\longrightarrow\Pi  \sigma_i^*(V) 
\buildrel\Pi\sigma_i^*(s_i)\over\longrightarrow\Pi^2 \sigma_i^*\sigma_i^*V 
\buildrel\simeq \over\longrightarrow \Pi^2 V =  V$$
is equal to $(-\id_V)$. 

\vskip .1cm

\noindent (d) For any $i=1, ..., n-2$ the two compositions
$$V\buildrel s_i\over\longrightarrow\Pi \sigma_i^*V
\buildrel\Pi\sigma_i^*(s_{i+1})\over\longrightarrow\Pi^2 \sigma_i^*\sigma_{i+1}^* V
\buildrel\Pi^2 \sigma_i^*\sigma_{i+1}^*(s_i) \over\longrightarrow\Pi^3 \sigma_i^*\sigma_{i+1}^*\sigma_i^* V,$$

$$V\buildrel  s_{i+1}\over\longrightarrow\Pi \sigma_{i+1}^*V
\buildrel\Pi\sigma_{i+1}^*(s_{i})\over\longrightarrow \Pi^2\sigma_{i+1}^*\sigma_{i}^* V
\buildrel\Pi^2 \sigma_{i+1}^*\sigma_{i}^*(s_{i+1}) \over\longrightarrow\Pi^3\sigma_{i+1}^*\sigma_{i}^*\sigma_{i+1}^* V$$
become equal after identifying the rightmost terms with
$$\Pi (\sigma_i\sigma_{i+1}\sigma_i)^*V = \Pi (\sigma_{i+1}\sigma_i\sigma_{i+1})^*V.$$

\vskip .1cm

\noindent (e) For any $i, j$ such that $|i-j|\geq 2$ the two compositions
$$V\buildrel s_i\over\longrightarrow\Pi  \sigma_i^*V \buildrel\Pi \sigma_i^*(s_j)\over\longrightarrow
\Pi^2 \sigma_i^*\sigma_j^*V, \quad \quad
V\buildrel s_j\over\longrightarrow\Pi \sigma_j^*V \buildrel\Pi \sigma_j^*(s_i)\over\longrightarrow
\Pi^2\sigma_j^*\sigma_i^*V $$
differ by sign, after identifying the rightmost terms with
$$(\sigma_i\sigma_j)^*V = (\sigma_j\sigma_i)^*V.$$

\end{271}
\vskip .1cm

\noindent Morphisms  are morphisms $V\to V'$ in $V^{\boxtimes_s n}$ commuting with the $s_i$.

\vskip .2cm

Thus the difference with $\Lambda^n_{\on{naive}}$ is that the $s_i$ are assumed to
be odd morphisms. 

\vskip .2cm

\noindent {\bf (3.6.2) Examples.} (a)
We have
$$\Lambda^n(\on{SVect^f}) \quad = \quad \on{Smod^f-}\kk[S_n]^c.$$
In particular,
$$\on{rk}\,\, K(\Lambda^n(\on{SVect^f})) \quad = \quad s(n).$$

\vskip .2cm
\noindent (b) 
Let $A$ be a finite-dimensional
$\kk$-superalgebra. We have then the crossed product superalgebra
$A^{\otimes n}\circledast \kk[S_n]^c$ generated by $A^{\otimes n}$ and odd 
generators  $s_i$, $i=1, ..., n-1$ which satisfy the
relations (2.4.1) together with the commutation relations
$$s_i\cdot a = a\cdot \sigma_i(a), \quad a\in A^{\otimes n}.$$
We have then
 $$\widehat\Lambda^n\(\on{Smod^f-}A\) \quad \simeq \quad 
\on{Smod^f-}\(A^{\otimes n}\circledast \kk[S_n]^c\).$$
\vskip 1cm

\vfill\eject

\section*{4. The categorical Koszul complexes}

\vskip .5cm

\noindent {\bf (4.1) Statement of the result.} Let $\mV$ be a superlinear category.
In this section we assume one of the following:

\vskip .2cm

(a) $\mV$ is abelian, and $\mV^{\widehat{\boxtimes}_sn}$ exists for every $n$.

\vskip .1cm

(b) $\mV$ is a 2-periodic perfect dg-category.

\vskip .2cm

Under these assumptions, we have the completed symmetric and exterior powers
$\widehat{\on{Sym}}^n\mV$, $\widehat{\Lambda}^n\medspace\!\mV$, as well as their products
$$\widehat{\Lambda}^p\medspace\!\mV\medspace\widehat\boxtimes_s\medspace \widehat{\Sym}^q\mV, \quad p,q\geq 0. 
\leqno (4.1.1)$$
The products in (4.1.1) are defined directly as consisting of
objects of $V^{\widehat\boxtimes_s (p+q)}$ with two types of equivariance
data commuting with each other: straight equivariance (2.2.1)
with respect to $S_q\subset S_{p+q}$ permuting the last $q$ 
factors, and ``twisted'' equivariance (3.6.1) with respect to
$S_p\subset S_{p+q}$ permuting the first $p$ factors. These products
are again abelian superlinear in the case (a) and 2-periodic perfect
in the case (b), so their Grothendieck groups are defined.
\vskip .1cm

\newtheorem*{312}{(4.1.2) Theorem}
\begin{312}
(a) There exist sequences of functors (categorical Koszul complexes)
\medskip

\begin{center}
\begin{tikzpicture}                                                              
\matrix(m)[matrix of math nodes, text height=1.5ex, row sep=1.5em, 
%column sep=-1em,
text depth=0.25ex]  
{
0 &{\phantom{X}}&  \widehat{\Lambda}^n\medspace\!\mV &{\phantom{iX}} &
(\widehat{\Lambda}^{n-1}\medspace\!\mV)
\medspace\widehat\boxtimes_s \medspace\mV& {\phantom{XX}}&
(\widehat{\Lambda}^{n-2}\medspace\!\mV)\medspace \widehat\boxtimes_s\medspace
(\widehat{\Sym}^2\mV)& {\phantom{iX}} & ...& {\phantom{X}}
& \widehat{\Sym}^n\mV&\phantom{X}& 0\\};

\draw[->,font=\scriptsize,>=angle 90] (m-1-1) -- (m-1-3);     
\draw[->,font=\scriptsize,>=angle 90] (m-1-3) -- node[auto]                     
{$D$} (m-1-5);                                     
\draw[->,font=\scriptsize,>=angle 90] (m-1-5) -- node[auto] {$D$}
(m-1-7);       
\draw[->,font=\scriptsize,>=angle 90] (m-1-7) -- node[auto] {$D$}
(m-1-9);        
\draw[->,font=\scriptsize,>=angle 90] (m-1-9) -- node[auto] {$D$}
(m-1-11);        
\draw[->,font=\scriptsize,>=angle 90] (m-1-11) -- 
(m-1-13);        

\end{tikzpicture}                                                                
  
\end{center}
\medskip
\begin{center}
\begin{tikzpicture}                                                              
\matrix(m)[matrix of math nodes, text height=1.5ex, row sep=2em, 
%column sep=-1em,
text depth=0.25ex]  
{
0&\phantom{X}&\widehat{\Sym}^n\mV&\phantom{iX}&
(\widehat{\Sym}^{n-1} \mV)
\medspace\widehat\boxtimes_s\medspace \mV&\phantom{XX}&
(\widehat{\Sym}^{n-2} \mV)\medspace\widehat\boxtimes_s \medspace
(\widehat{\Lambda}^2\medspace\!\mV) &\phantom{X}& ... &\phantom{X}& 
\widehat{\Lambda}^n\medspace\!\mV&\phantom{X}& 0\\};

\draw[->,font=\scriptsize,>=angle 90] (m-1-1) -- (m-1-3);     
\draw[->,font=\scriptsize,>=angle 90] (m-1-3) -- node[auto]                     
{$\Delta$} (m-1-5);                                     
\draw[->,font=\scriptsize,>=angle 90] (m-1-5) -- node[auto] {$\Delta$}
(m-1-7);       
\draw[->,font=\scriptsize,>=angle 90] (m-1-7) -- node[auto] {$\Delta$}
(m-1-9);        
\draw[->,font=\scriptsize,>=angle 90] (m-1-9) -- node[auto] {$\Delta$}
(m-1-11);        
\draw[->,font=\scriptsize,>=angle 90] (m-1-11) -- 
(m-1-13);        
\end{tikzpicture}                                                                
  
\end{center}

(b) These functors are exact, if $\mV$ is abelian, and pre-exact, if $\mV$
is perfect. On the level of Grothendieck groups, these functors give rise to
chain complexes:
$$D_*^2=0,\quad  \Delta_*^2=0.$$

(c) The induced complexes of complexified Grothendieck groups are acyclic. 
\end{312}

\vskip .2cm

%If all these maps are isomorphisms, then the complexes of the Grothendieck
%groups from (3.1.2)(b) look even more like Koszul complexes. 
%
%\vskip .2cm

If the product maps
$$
  \longmap{\circledast}
  {K_\CC^\bullet(\widehat{\Lambda}^p\medspace\! \mV)\otimes
  K^\bullet_\CC(\widehat{\Sym}^q \mV)}
  {K^\bullet_\CC(\widehat{\Lambda}^p\medspace\!\mV\medspace\widehat\boxtimes_s\medspace 
  \widehat{\Sym}^q\mV)}
$$ 
of (1.7.3) are isomorphisms,  then the complexes of the Grothendieck
groups in (4.1.2) look even more like Koszul complexes.

\vskip .2cm
\noindent {\bf (4.1.3) Example.} Let $A$ be a finite-dimensional $\kk$-superalgebra,
and $\mV$ be the abelian category $\on{Smod^f-}A$. Then
 $$\widehat{\Sym}^n(\mV) = \on{Smod^f-}(A^{\otimes n}[S_n]), \quad 
\quad  \widehat{\Lambda}^n(\mV) = \on{Smod^f-}(A^{\otimes n}\circledast \kk[\Sn]^c),$$
see Examples (1.6.7), (2.2.8)(b) and (3.6.2)(b). In this case the
homomorphisms (1.7.3) are isomorphisms (see Example (1.7.4)(c)),
and we get the equality 
%\marginpar{Wrong? (but fin.dim $\Rightarrow$ simple is true).}
$$\biggl(\sum_{n\geq 0} \dim \, \, K^\bullet_\CC(\on{Smod^f-}A^{\otimes
  n}[S_n]) q^n\biggr) \cdot 
\biggl( \sum_{n\geq 0} \dim \, \, K^\bullet_\CC(\on{Smod^f-}(A^{\otimes
  n}\circledast \kk[\Sn]^c))(-q)^n\biggr) =1.$$

\vskip .3cm

\noindent {\bf (4.2) Partial symmetrization.} Let $G$ be a finite group
acting on a linear category $\mW$, so we have functors $\varrho(g): \mW\to \mW$
as in (2.1). Let $H$ be a subgroup of $G$.
Then the forgetful functor from the category of $G$-equivariant
to that of $H$-equivariant objects has a left adjoint, 
\begin{eqnarray*}
I_H^G:  \mW^H&\longrightarrow& \mW^G  \\
W&\longmapsto& \bigoplus_{[g]\in G/H} \varrho(g)(W),
\end{eqnarray*}
called
{\em `partial symmetrization'}.
Here $g$ runs over representative classes of right cosets of $G$ by
$H$, and the
direct sum carries the obvious $G$-equivariant structure. 
\vskip .2cm

\noindent {\bf (4.2.1) Examples.} (a)  Suppose that the $G$-action on $\mW$ is trivial.
Then $W$ is a representation of $H$ in $\mW$, and $I_H^G(W) = \on{Ind}_H^G(W)$
is the induced representation. 

\vskip .1cm

\noindent
(b) Let $G= S_{q+1}$ be the group of permutations of $\{0,1, ..., q\}$, and $H=S_q$
be the group of permutations of $\{1, ..., q\}$. Then a set of representatives
for $S_{q+1}/S_q$ is provided by the transpositions
$$\id, \, (0,1), \,\, (0,2), \,\, ..., \,\, (0,q),$$
so for an $S_q$-equivariant object $W$ its partial symmetrization is given by
$$I_{S_q}^{S_{q+1}}(W) \quad =\quad W\,\, \oplus \,\,\bigoplus_{i=1}^q \varrho((0,i))(W).$$

\vskip .1cm
\noindent
(c) Let $H=1$ is the trivial group. The functor $I_1^G$ plays a
crucial role in the proof of the main theorem in
\cite{innerproducts}, where $I_1^G$ is denoted $A'$. 

\vskip .2cm 
The following is obvious from the definition of $I_H^G(W)$ as a direct sum.

\newtheorem*{322}{(4.2.2) Proposition}
\begin{322}
If $\mW$ is abelian, then $I_H^G$ is exact. If $\mW$ is a perfect dg-category,
then $I_H^G$ is pre-exact.
\end{322}

\vskip .3cm

\noindent {\bf (4.3) The functors $D_{p,q}$ and $\Delta_{p,q}$.}
Denote by
$$J_p: \widehat{\Lambda}^p\medspace\!\mV\longrightarrow
(\widehat{\Lambda}^{p-1}\medspace\!\mV)\medspace\widehat\boxtimes_s\medspace \mV
$$
the functor of forgetting the part of the twisted equivariance data not
pertaining to $S_{p-1}$. In other words, an object of $\widehat{\Lambda}^{p}(\mV)$
is an object $V\in \mV^{\widehat\boxtimes_s p}$ with isomorphisms $s_i: V\to \Pi\sigma_i^*(V)$,
$i=1, ..., p-1$, satisfying the conditions in (3.6.1). By considering only
$s_1, ..., s_{p-2}$, we get an object of $\widehat{\Lambda}^{p-1}(\mV)\widehat\boxtimes_s \mV$,
which we denote $J_p(V)$. We now define the functor $D_{p,q}$ to be the composition
\noindent
\begin{center}
\begin{tikzpicture}                                                              
\matrix(m)[matrix of math nodes, column sep=2cm, text depth=0.25ex]  
{
\widehat{\Lambda}^p\medspace\!\mV\medspace\widehat\boxtimes_s
\medspace\widehat{\Sym}^q\mV
&
(\widehat{\Lambda}^{p-1}\medspace\!\mV)\medspace
\widehat\boxtimes_s\medspace \mV \medspace \widehat\boxtimes_s\medspace
(\widehat{\Sym}^q\mV)
&
(\widehat{\Lambda}^{p-1}\medspace\!\mV)\medspace\widehat\boxtimes_s\medspace
(\widehat{\Sym}^{q+1}\mV).\\
};

\draw[->,font=\scriptsize,>=angle 90] (m-1-1) -- node[auto]
{$J_p\widehat\boxtimes_s\id$} (m-1-2);     
\draw[->,font=\scriptsize,>=angle 90] (m-1-2) -- node[auto]                     
{$\id\widehat\boxtimes_s I_{S_q}^{S_{q+1}}$} (m-1-3);           
\end{tikzpicture}                                                                
  
\end{center}

Similarly, we have the functor
$$\widetilde{J}_q: \widehat{\Sym}^q\mV\longrightarrow
\mV\medspace\widehat\boxtimes_s \medspace
\widehat{\Sym}^{q-1}\mV
$$ 
 forgetting the part of equivariance structure not
pertaining to $S_{q-1}$. Again, we have the partial symmetrization functor
$${\widetilde{I}}_{S_{p-1}}^{S_p}: \widehat{\Lambda}^p(\mV)\medspace \widehat\boxtimes_s\medspace\mV = 
\bigl(\mV^{\widehat\boxtimes_s (p+1)}\boxtimes_s \on{Sgn}\bigr)^{S_p} \to 
\bigl(\mV^{\widehat\boxtimes_s (p+1)}\boxtimes_s \on{Sgn}\bigr)^{S_{p+1}} = \widehat{\Lambda}^{p+1}(\mV).
$$
Explicitly, for an object $W$ of $\widehat{\Lambda}^p(\mV)\widehat\boxtimes_s\mV$ we have
$${\widetilde{I}}_{S_p}^{S_{p+1}}(W) \quad =\quad W\,\,\oplus \,\, \bigoplus_{i=1}^p \, \Pi \varrho((i, p+1))(W),
$$
where $\rho$ is the permutation 2-representation.
We define the functor $\Delta_{p,q}$ to be the composition
\begin{center}
\begin{tikzpicture}                                                              
\matrix(m)[matrix of math nodes, column sep=2cm, text depth=0.25ex]  
{
\widehat{\Lambda}^p\mV\medspace\widehat\boxtimes_s\medspace
\widehat{\Sym}^q\mV
&
(\widehat{\Lambda}^p\medspace\!\mV)\medspace\widehat\boxtimes_s\medspace
\mV\medspace \widehat\boxtimes_s\medspace(\widehat{\Sym}^{q-1}\mV)
&
(\widehat{\Lambda}^{p+1}\medspace\!\mV)\medspace\widehat\boxtimes_s
\medspace(\widehat{\Sym}^{q-1}\mV).  
\\
};

\draw[->,font=\scriptsize,>=angle 90] (m-1-1) -- node[auto]
{$\widetilde{J}_q\widehat\boxtimes_s\id$} (m-1-2);     
\draw[->,font=\scriptsize,>=angle 90] (m-1-2) -- node[auto]                     
{$\id\widehat\boxtimes_s \widetilde{I}_{S_{p-1}}^{S_p}$} (m-1-3);           
\end{tikzpicture}                                                                
  
\end{center}

It is clear that the $D_{p,q}$ and $\Delta_{p,q}$ are exact if $\mV$ is abelian
and pre-exact if $\mV$ is perfect. So they induce morphisms
$(D_{p,q})_*$ and $(\Delta_{p,q})_*$ of the Grothendieck groups.

\vskip .2cm

\newtheorem*{331}{(4.3.1) Proposition}
\begin{331}
We have the identities:
\begin{eqnarray*}
(D_{p-1, q+1})_* (D_{p,q})_*&=&0,\\
(\Delta_{p-1, q+1})_* (\Delta_{p,q})_*&=&0, \quad\text{and}\\
(\Delta_{p-1, q+1})_* (D_{p,q})_* + (D_{p+1, q-1})_* (\Delta_{p,q})_*
&=& (p+q)\cdot \id.
\end{eqnarray*}
\end{331}

Thus the first two series of equalities mean that we have chain complexes,x while the
third series implies that the complexes become acyclic
after tensoring with $\CC$. Indeed, they mean that $\Delta_*$ is 
a contracting homotopy for $D_*$ and vice versa: the isomorphism of the $n$th complex
given by multiplication by $n=p+q$ is homotopic to zero. Thus Proposition 4.3.1
implies Theorem 4.1.2. 

\vskip .2cm

\noindent {\sl Proof of (4.3.1):} We start with the first identity. It
is enough to 
consider $p=2$, the first $p-2$ arguments being immaterial. 
The composition
\begin{center}
  \begin{tikzpicture}
    \matrix(m)[matrix of math nodes, column sep=2cm]
    {
    (\widehat{\Lambda}^2\mV)\medspace \widehat\boxtimes_s \medspace
    (\widehat{\Sym}^q\mV) 
    &
    \mV\medspace\widehat\boxtimes_s\medspace(\widehat{\Sym}^{q+1}\mV)
    &
    \widehat{\Sym}^{q+2}\mV
    \\
    };

    \draw[->,font=\scriptsize,>=angle 90] (m-1-1) -- node[auto] 
    {$D_{2,q}$} (m-1-2);
    \draw[->,font=\scriptsize,>=angle 90] (m-1-2) -- node[auto]{$D_{1,q+1}$}
    (m-1-3); 
  \end{tikzpicture}
\end{center}
is just the partial symmetrization $I_{S_q}^{S_{q+2}}$.
Here $S_q$ acts on $\{1, ..., q+2\}$ by fixing the first two
elements. Hence the elements of a coset $\sigma S_q$ are exactly those
elements $\tau\in S_{q+2}$ with $\tau(1)=\sigma(1)$ and $\tau(2)=\sigma(2)$.
So, any system of coset-representatives $\{\sigma_{i,j}\}$ may be labeled by
ordered pairs $i,j$ in $\{1,\dots,q+2\}$, 
with $\sigma_{i,j}(1)=i$ and $\sigma_{i,j}(2)=j$. Note that we can
choose these representatives in
such a way that we have
$$\sigma_{i,j}=\sigma_{j,i}\circ(1,2).$$
Let  $V$ be an object 
 of $(\widehat{\Lambda}^2\medspace\!\mV)\medspace\widehat\boxtimes_s\medspace(
 \widehat{\Sym}^2\mV)$. Then $V$  
is twisted equivariant in first two factors, and we have
$$\varrho((1,2))(V) \quad \simeq\quad \Pi V,$$
and more generally,
$$
  \varrho(\sigma_{i,j})(V)\cong\Pi\varrho(\sigma_{j,i}(V)).
$$
Hence their classes in the $K$-group cancel, and 
$$(D_{1,q+1})_*(D_{2,q})_* \langle V\rangle =0.$$
The second identity, 
$(\Delta_{1,q+1})_*(\Delta_{2,q})_* =0,$
is proved in a similar way.

\vskip .2cm

We now prove the third identity. Let $V$ be an object of
$\widehat{\Lambda}^p(\mV)\medspace\widehat\boxtimes_s\medspace
\widehat{\Sym}^q(\mV)$. 
As in (4.1), we view $V$ as on object of $\mV^{\widehat\boxtimes_s (p+q)}$ with
two types of equivariance data. Then
$$
  \Delta_{p-1, q+1} (D_{p,q}(V)) \cong
  \(V\oplus\bigoplus_{i=1}^{p-1}\Pi(\varrho(i,p)V)\)
  \oplus
  \bigoplus_{j=1}^q\(\varrho(p,p+j)V\oplus\bigoplus_{i=1}^{p-1}
  \Pi(\varrho(i,p,p+j)V) \).
$$ 
By virtue of the two kinds of equivariance data of $V$,
the first $p$ summands each are isomorphic to $V$, each of the next $q$
summands is isomorphic to one of the form $\varrho(p,p+j,p+1)$, and
the remaining $q(p-1)$ summands are isomorphic to ones of the form
$\varrho(i,p+j,p+1)$. Here the last step used that
$$(i,p,p+j)=(i,p+j,p+1)\circ(p+1,p+j)(i,p).$$  
Likewise, $D_{p+1, q-1}(\Delta_{p,q}(V))$ is a direct sum
of $(p+1)q$ terms, $q$ of which are isomorphic to $V$. The remaining
$pq$ terms are $\Pi(\varrho(i,p+j,p+1)$ with $1\leq i\leq
p$ and $1\leq j\leq q$, pairing up with the corresponding terms in $
\Delta_{p-1, q+1} D_{p,q}(V)$ to cancel in the Grothendieck group.
\qed

\vskip 1cm

\vfill\eject

\section*{5. Examples from representation theory}

\vskip .5cm
Let $\mV$ be a superlinear category of one of the types considered in
(4.1). The spaces 
$$\bigoplus_{n\geq 0} \, K^\bullet_\CC(\widehat{Sym}^n\mV) \quad\text{
  and }\quad \bigoplus_{n\geq 0}
\, K^\bullet_\CC (\widehat{\Lambda}^n\mV)\leqno (5.1)$$
can be thought as  analogues of the bosonic and fermionic Fock spaces. 
Several recent papers realize these spaces as basic representations of
some infinite-dimensional Lie algebras. In each of these situation the first space
corresponds to the ``untwisted" case, while the second one corresponds 
to an appropriately   ``twisted' case. In this section we recall some examples
and give their interpretation from our point of view.

\vskip .3cm

\noindent {\bf (5.2) Wedge and spin representations.}
(a) Take $\mV= \on{SVect^f}$ to be the category of finite-dimensional
super $\kk$-vector spaces. Then ${\Sym}^n\mV = \on{Rep}_\kk(S_n)$
is the category of finite-dimensional $\ZZ/2$-graded 
representations of $S_n$ over $\kk$. 
Here it is not necessary to complete symmetric and exterior
powers. There are no self-associate irreducible objects in
$\on{Rep}_\kk(S_n)$, and  
$$K^\bullet_\CC({\Sym}^n\mV)=K^0_\CC(\Sym^n\mV)$$ 
is the complexified representation ring of
$S_n$. Denote this ring $\mathcal{R}_\kk(S_n)$. 

Assume $\on{char}(\kk)=0$. Then a basis of 
$\mathcal{R}_\kk(S_n)$ is labelled by partitions
$\alpha = (\alpha_1 \geq ... \geq \alpha_l)$ of $n$
which label irreducible representations.
In fact, 
we have an identification
$$
  \bigoplus_{n\geq 0} \, \mathcal R_\kk(S_n) \!\!\quad\simeq\quad\!\! 
  \mathfrak F
$$
%\Lambda^{\infty/2} \,\, \CC\ls{t}$$
where the right-hand side is the basic
projective representation of the
Lie algebra $$\mathfrak{gl}(\infty) := \{(a_{i,j})_{i,j\in\ZZ} \mid
\text{almost all $a_{i,j}$ equal $0$}\},$$
see \cite[Lecture 4]{Kac:Raina} 
or \cite[Rem. 2.3.13]{Kleshchev}. Here $\mathfrak F$ stands for ``Fock''.
To describe $\mathfrak F$, let $V=\CC[t,t\inv]$ be the vector space of
all finite formal Laurent polynomials over $\CC$. 
Then $\mathfrak{gl}_\infty$ acts on $V$ in the standard way: let
$E_{i,j}$ be the matrix with only one entry, $1$ in the spot $(i,j)$. Then
$E_{i,j}t^k=\delta_{j,k}t^i$. The semi-infinite wedge
$\Lambda^{\infty/2}V$ inherits a $\mathfrak{gl}_\infty$-action from
$V$, and $\mathfrak F$ is the smallest subrepresentation containing
the {\em `vacuum vector'}
$$
  1\wedge t\inv\wedge t^{-2}\wedge t^{-3}\wedge\dots
$$
A basis of $\mathfrak F$
is given by semiinfinite wedge products of monomials in $t$:
$$ t^{i_0}\wedge t^{i_1} 
\wedge t^{i_2}\wedge ..., \quad i_0>i_1..., \quad i_k=-k, \,\, k \gg 0.$$
Now, given such a wedge, we form a partition
$$\alpha = (\alpha_0\geq \alpha_1 \geq ...), \quad \alpha_k = k+i_k,$$
and this gives a bijection between the bases of the two spaces.
In fact, the Young graph of partitions and their inclusions is
interpreted as the ``crystal graph'' of the wedge representation.
It means that the edges of the Young graph
(which correspond to adding/removing a node from a Young diagram)
 exactly describe  the action of the $E_{ij}$ on the basis wedge products).

\vskip .2cm

(b) Keeping the assumption $\on{char}(\kk)=0$, we have 
the following companion description for the exterior powers:
$$\bigoplus_n \, K^\bullet_\CC({\Lambda}^n(\mV))\!\!\quad \simeq \quad\!\!\Sigma,$$
where $\Sigma$ is the basic ``spin'' projective representation of the
Lie algebra $\mathfrak{o}(\infty)$. The direct sum on the
left is the direct sum of the Grothendieck groups of the categories
$\on{Smod}^f-\kk[\Sn]^c$ of projective  super-representations
of $S_n$. The representation $\Sigma$ was defined in
\cite{Kac:Peterson}, while its interpretation in terms of 
projective representations of symmetric groups 
was studied in \cite{Jing},
using the techniques of vertex operators. 

\vskip .3cm

\noindent{\bf (5.3) The Kac-Moody algebras $A_{p-1}^{(1)}$ and $A_{p-1}^{(2)}$.}
Consider now the same situation as in (5.2) except assume that 
$\on{char}(\kk)=p>0$. 

\vskip .2cm

Irreducible representations of $S_n$ over $\kk$ are now labelled
by $p$-regular partitions, i.e., partitions where no part
is repeated $p$ or more times, see \cite{Kleshchev}
for details. For example, 2-regular
means strict. The spaces (5.1) in this case are identified with
the fundamental (projective) representations of the Kac-Moody Lie algebras 
$$A_{p-1}^{(1)} = \mathfrak{sl}_p \, (\CC[t, t^{-1}]), \quad A_{p-1}^{(2)} = 
\bigl\{g(t)\in \mathfrak{sl}_p\, \CC[t, t^{-1}]: \, g(t^{-1}) = g(t)^*\bigr\}.
\leqno (5.3.1)$$
Here $*$ means transposition. 
The symmetric powers in (5.1) make sense for all $p$,
while the exterior powers are defined only for $p>2$, in the same way
as the Lie algebras in (5.3.1). 
See \cite{Kleshchev} for more details and references.

\vskip .3cm

\noindent {\bf (5.4) Twisted and untwisted Kac-Moody algebras
of type ADE.} Let $\on{char}(\kk)=0$ and let $\Gamma$ be a finite 
group isomorphic to a subgroup of $SL_2(\CC)$. By the McKay correspondence,
$\Gamma$ corresponds to a root system of type ADE. Denote
by $\mathfrak{g}_\Gamma$ the corresponding semisimple split
Lie algebra over $\CC$. We have then the Kac-Moody algebras
$\mathfrak{g}_\Gamma^{(1)}, \mathfrak{g}_\Gamma^{(2)}$
defined as in (5.3.1) but with the transposition replaced
by the Cartan involution.

On the other hand, take $\mV= \on{Rep}(\Gamma)$ to be
the category of finite-dimensional
representations of $\Gamma$. 
Then
$$\widehat{Sym}^n(\mV) = \on{Rep}(S_n\smallint \Gamma^n),
\quad \widehat{\Lambda}^n(\mV) = \on{Rep}_c(S_n\smallint\Gamma^n)$$
are the categories of the all representations (resp. projective
representations of a certain kind) of the wreath product
$S_n\smallint\Gamma^n$. In fact, since we work in characteristic
zero, all the categories involved are semisimple and we can drop
the hats over symmetric and exterior powers.

The Grothendieck groups of such representations
were studied in \cite{Frenkel:Jing:Wang}, 
\cite{Frenkel:Jing:Wang:Twisted}, whose results mean that
$$\bigoplus_{n\geq 0} \, K^\bullet_\CC(\Sym^n(\mV)) = \mathfrak{F}_\Gamma^{(1)}, \quad
\bigoplus_{n\geq 0} \, K_\CC^\bullet(\Lambda^n(\mV)) =
\mathfrak{F}_\Gamma^{(2)}$$ 
are the fundamental representations of $\mathfrak{g}_\Gamma^{(1)}$
and $\mathfrak{g}_\Gamma^{(2)}$ respectively.

\vskip 1cm

\vfill\eject

\section*{6. Effect on categorical characters}

\vskip .5cm
Let $\varrho$ be a linear 2-representation of $G$ on
$\mV=\on{Vect^m}$. By \cite{Kapranov:Voevodsky}, we may view
$\varrho$ as a matrix-representation: $\varrho(g)$ is described by
an $m\times m$-matrix $(V_{i,j})$ whose entries are vectorspaces, and
the $\phi_{g,h}$ and $\phi_1$ amount to matrices of isomrophisms
between the relevant entries. 

\vskip .2cm
{\bf (6.1) Symmetric powers.} The symmetric powers $\on{Sym^n}(\varrho)$ are again
2-representations of $G$, and their categorical
characters are described as follows: for $g\in G$, we have an
isomorphism of $\mathbb N$-graded $C_g$-representations
$$
  \bigoplus_{n\geq 0}X_{\on{Sym^n}(\varrho)}(g)q^n
  \medspace\cong\medspace
  \bigotimes_{k\geq 1} S_{q^k}\(X_{\varrho}(g^k)^{\ZZ/k\ZZ}\). \leqno (6.1.1)
$$
Here $q$ is a dummy variable, and $$S_q(V)=\bigoplus \on{sym^m}(V)q^m$$ is
the total symmetric power of the vector space $V$. 
\begin{proof}[Proof of (6.1.1)]
  The $n$th tensor power $\mV^{\boxtimes_sn}$ carries an action of
  $S_n\smallint G$. The group $G$ still acts diagonally on the category
  $\on{Sym^n}(\mV)$ of $S_n$-equivariant objects.
  By \cite[Theorem 5.13]{innerproducts}, we have an isomorphism of vector
  spaces 
  $$
    X_{\on{Sym^n}(\varrho)}(g) \cong  \(\bigoplus_{\sigma\in
      S_n}X_{\varrho^{\boxtimes n}}(\sigma;g,\dots,g)\)^{S_n}.
  $$
  Fix $\sigma\in S_n$ and let $n_k$ be the number of $k$-cycles in
  $\sigma$. 
  Then the action of $(\sigma;g,\dots,g)$ on $\mV^{\boxtimes
  n}$ is isomorphic to the (external) categorical tensor product of
  auto-equivalences $(c_k;g,\dots,g)$ of $\mV^k$, where $c_k$
  acts as a cyclic permutation of the tensor factors.
  Since the categorical trace is multiplicative (see
  \cite[Thm. 2.5]{innerproducts}), we are reduced to
  calculating the trace of $(c_k;g,\dots,g)$.
  An examination of
  the diagonal entries of the relevant matrix yields an isomorphism 
  $$X_{\varrho^{\boxtimes k}}(c_k;g,\dots,g)\cong X_{\varrho}(g^k).$$
  More precisely, both sides are isomorphic to 
  $$
    \bigoplus_{\ul i}V_{i_1,i_2}\tensor V_{i_2,i_3}\tensor\dots\tensor
    V_{i_k,i_1}, 
  $$
  and $\langle
  c_k\rangle$ acts by cyclic permutation of the tensor factors.
%  Here the left-hand side is a trace in $(\SVect)^{m^k}$, while the
%  right-hand side is a trace in $(\SVect)^m$.
  Taking fixed points under the action of the centralizer 
  $$C_\sigma\cong\prod_k S_{n_k}\smallint (\ZZ/k\ZZ),$$ we obtain
  \begin{eqnarray*}
    X_{\varrho^{\boxtimes n}}(\sigma;g,\dots,g)^{C_\sigma}&% \cong&
                                % \bigotimes_{k} 
%    \(\(\ttr(g^k)^{\ZZ/k\ZZ}\)^{\tensor n_k}\)^{S_{n_k}}\\
   \cong& \bigotimes_{k}\on{sym^{n_k}}\(X_{\varrho}(g^k)^{\ZZ/k\ZZ}\), 
  \end{eqnarray*}
  and hence 
  $$
    X_{\on{Sym^n}(\varrho)}(g)\cong \bigoplus_{n=\sum kn_k}
    \bigotimes_{k\geq 1}\on{sym^{n_k}}\(X_{\varrho}(g^k)^{\ZZ/k\ZZ}\).
  $$
  This is the coefficient of $q^n$ on the right-hand side of (6.1.1).
\end{proof}

\vskip .2cm
{\bf (6.2) Exterior powers.} 
Let now $\varrho$ be a superlinear 2-representation of $G$ on
$\mV=\on{(\SVect)^m}$. Just as in the linear case, we may view
$\varrho$ as a matrix-representation: $\varrho(g)$ is described by
an $m\times m$-matrix of supervectorspaces, and the $\phi_{g,h}$ and
$\phi_1$ amount to even maps between such matrices. The formalism
of \cite{Kapranov:Voevodsky} goes through in the super setting. For
instance, composition of superlinear functors corresponds to matrix multiplication,
with $+$ and $\cdot$ replaced by $\oplus$ and $\tensor$ of
supervectorspaces. 

To calculate the categorical supercharacter of the total exterior
power of $\varrho$, we apply the same argument as in the proof of
(6.1.1), but this 
time with a different $S_n$-action. This yields the formula
$$
  \bigoplus_{n\geq 0}\on{sX^\bullet_{\Lambda^n(\varrho)}}(g)q^n\cong
  \bigotimes_{k\geq1}S_{q^{2k-1}}\(\on{sX^\bullet_\varrho}(g^{2k-1})^+\)\tensor
  \Lambda_{q^{2k}}\(\on{sX^{\bullet+\bar1}_\varrho}(g^{2k})^-\).
  \leqno{(6.2.1)}
$$
For odd $k$, the cyclic group $\ZZ/k\ZZ$ acts 
on $sX^\bullet(g^k)$ as in the proof of (6.1.1). 
For even $k$, the action acquires the sign
$\chi_{\on{spin}}$ (see Definition (3.4.8)). The superscripts $+$ and
$-$ indicate the invariant parts under these respective
actions, and $\Lambda_q$ denotes the total exterior power of
supervectorspaces. 

\vfill\eject

\section*{7. Open questions}

\noindent {\bf (7.1) Symmetric
  product orbifolds with discrete torsion.}
With the corrected $\epsilon$ of Section (3.4.7),  Dijkgraaf's Formula
\cite[(4.7)]{Dijkgraaf}
for the Hilbert space of the
symmetric product orbifold with discrete torsion
becomes  
%I don't entirely understand this part of the paper, but believe that
%there might be an additional typo in \cite[(4.7)]{Dijkgraaf}, since
%the numbers $N_n$ do not turn up in the tensor factors on the
%right-hand side. Looking at the argument in (3.29) of
%\cite{DMVV}, my guess is that 
%in \cite[(4.7)]{Dijkgraaf}, 
%$S_{p^{2n-1}}$ should be $S^{N_{2n-1}}$ and $\Lambda_{p^{2n}}$ should
%be $\Lambda^{N_{2n}}$. 
%With the corrected $\epsilon$,
%The right-hand side of \cite[(4.7)]{Dijkgraaf} becomes
$$
  \mathcal H^c(S^NX) = \bigoplus_{\stackrel{\text{even }\{N_n\}}{N=\sum nN_n}}
  \bigotimes_{n>0}\left(
  S^{N_{2n-1}}\mathcal H_{2n-1}^+\otimes  \Lambda^{N_{2n}}\mathcal
  H_{2n}^- 
  \right)\oplus
  \bigoplus_{\stackrel{\text{odd }\{N_n\}}{N=\sum nN_n}}
  \bigotimes_{n>0}
  \Lambda^{N_{n}}\mathcal H_{n}^+.
  \leqno{(7.1.1)}
$$
In other words, if the overall parity of $\sigma$ is odd then long
strings of odd length exhibit fermionic, not bosonic behaviour.
There is no obvious way to summarize (7.1.1) into a 
product formula for all $N$, as in \cite{DMVV}.
One would expect this to change when the super point of view is taken
into account: assume that Dijkgraaf's Hilbert space $\mathcal H$ is
$\ZZ/2$-graded by the Fermion number. 
\newtheorem*{612}{(7.1.2) Question}
\begin{612}
  Can one incorporate the
  $\ZZ/2$-grading of $S_n$ and the Koszul sign
  $$(-1)^{\on{deg}(\sigma)\on{deg}(\tau)}$$ into the CFT-interpretation
  \cite[4.2]{Dijkgraaf} such that (6.1.1)
  becomes
  $$
    \mathcal H_\bullet^c(S^NX) = \bigoplus_{N=\sum nN_n}
    \bigotimes_{n>0}
    S^{N_{2n-1}}\(\mathcal H_{\bullet,2n-1}^+\)\otimes
    \Lambda^{N_{2n}}\(\mathcal H_{\bullet+1,2n}^-\)\quad?
    \leqno{(7.1.3)}
  $$
  Here $\mathcal H^c_\bullet\(S^NX\)$ is the super-Hilbertspace of the
  symmetric product orbifold with discrete torsion.%, and
%  $\mathcal H_{\bullet,2n}^-$ is to be adequately reinterpreted. 
\end{612}
Formula (7.1.3) is obtained by replacing $\epsilon$
with $e$ in Dijkgraaf's argument.
Note that with this adjustment, $n$-cycles behave as
``bosons'' if $n$ 
is odd and as ``fermions'' if $n$ is even, as postulated in
\cite[p.7]{Dijkgraaf}, but the distinction between even and odd
partitions in \cite[(4.7)]{Dijkgraaf} has vanished. 

Further, one would expect a
Dijkgraaf-Moore-Verlinde-Verlinde type product formula for the
super-Hilbertspace of the total symmetric product orbifold with
discrete torsion:
$$
  \sum_{n\geq_0}\mathcal H^c_\bullet\(S^n\)q^n = 
  \bigotimes_{k\geq 1}S_{q^{2k-1}}\(\mathcal H_{\bullet,2k-1}^+\)
  \tensor
  \Lambda_{q^{2k}}\(\mathcal H_{\bullet+1,2n}^-\).
$$
Here $q$ is a formal variable and $S_q$ and $\Lambda_q$ denote the total
symmetric and exterior power of super vectorspaces.

\vskip .3cm

\noindent {\bf (7.2) Operations on the Bondal-Larsen-Lunts ring.}
Our results suggest that the Grothendieck ring $\mathcal{PT}$ of
pretriangulated categories defined in \cite{Bondal:Larsen:Lunts}
has well-defined nonlinear operations $\Sym^n$.
 The relations in $\mathcal{PT}$
come from semi-orthogonal decomposition, so a study of these
operations would need to address the following question.
\newtheorem*{621}{(7.2.1) Question}
\begin{621}
  How do the $\Sym^n$ interact with semi-orthogonal decompositions?
\end{621}
\newtheorem*{622}{(7.2.2) Question}
\begin{622}
  Is there a way to incorporate the super point of view and exterior
  powers into this picture?
\end{622}
It is natural to expect that the formalism for these operations is
similar to that of a 2-special lambda ring, see \cite{2lambdarings} or
\cite{Rezk}. 
%\newtheorem*{51}{(5.1) Proposition}
%\begin{51}
%Let $\mV$ be a pretriangulated category, and suppose that
%$\mV', \mV''$ form a semi-orthogonal decomposition
%of $\mV$. Then $\widehat{\Sym}^n(\mV)$
%has a semiorthogonal decomposition by categories
%quasi-equivalent to
%$$\widehat{\Sym}^n(\mV'), \quad  \widehat{\Sym}^{n-1}(\mV')\
%\widehat\boxtimes_s \mV'', \quad  \widehat{\Sym}^{n-2}(\mV')\widehat\boxtimes_s
%\widehat{\Sym}^2(\mV'') , ...$$
%Similarly for $\widehat{\Lambda}^n(\mV)$. 
%\end{51}
%
%\vskip .2cm
Given a group $G$, one can consider pretriangulated
categories with $G$-action, so we get a kind of
2-representation ring $\mathcal{PT}_G$,
which should also have this kind of structure. Let $c_k$
be a cyclic permutation of length $k$, and note that the categorical
trace of the action of $(c_k;g,\dots,g)$ on $X^k$ involves the fixed
point set $\(X^k\)^{(c_k;g,\dots,g)}$. This can be identified with the
fixed points $X^{g^k}\sub X$, who are in turn involved in the
categorical trace $\ttr^\bullet(g^k)$. This observation motivates the
following question.
\newtheorem*{623}{(7.2.3) Question}
\begin{623}
  Is it possible to generalize the discussion of Section 6 to a wider
  class of 2-representations?
\end{623}
%Exactly what kind of lambda-structure we would
%figure out in the other paper devoted to
%2-special lambda rings.
%
%
%\marginpar{Question 16}

\vskip .3cm

\noindent {\bf (7.3) Elliptic genera.}
The expression $\bigotimes_{k\geq 1}S_{q^k}(-)$ in Formula (6.1.1) is
familiar from the definition of the Witten genus. Formula
(6.2.1) reminds of the ``expansion at the other cusp of the
equivariant signature of the loop space'' in
\cite[p.82]{Hirzebruch:Berger:Jung}.  
\newtheorem*{631}{(7.3.1) Question}
\begin{631}
  Can one use our constructions to define the Witten genus of $X$ in
  terms of $D^b(X)$, the bounded   derived category  
  of coherent sheaves on $X$ (or its pre-triangulated enhancement)?
\end{631}
Note that if we use the tensor structure on the derived category, then
an affirmative answer to this question follows indirectly from results
of Balmer 
\cite{Balmer}, who is able to recover the whole scheme $X$ from 
$D^b(X)$ as a tensor triangulated category. However, there are examples 
of pairs $X, X'$ of algebraic varieties (related by a birational ``flop") for which we have 
an equivalence $D^b(X)\simeq D^b(X')$ (not preserving tensor structures), 
as well as  an equality of the Witten genera. See \cite{Bridgeland, Totaro}.
This suggest that there may be a  relation between the (pre-)triangulated
category structure on  $D^b(X)$  and the
Witten genus.

\vskip 2cm

\vfill\eject

\end{document}